\documentclass[a4paper,10pt]{article}
\usepackage{amsmath}
\usepackage{amsthm}
\usepackage{amssymb}
\usepackage{mathrsfs}
\usepackage{subfigure}
\usepackage{fancyhdr}
\usepackage{latexsym}
\usepackage{graphicx}
\newcommand{\dimostr}{\noindent {\it Proof}\par\noindent}
\newtheorem{teorema}{Theorem}
\newtheorem{lemma}{Lemma}[section]
\newtheorem{proposizione}[lemma]{Proposition}
\newtheorem{corollario}{Corollary}[section]
\newtheorem{definizione}{Definition}[section]
\newtheorem{osservazione}{Remark}[section]

\title{On the regularization of the collision solutions of the one-center problem with weak forces\footnote{Work partially supported by PRIN Project \textquotedblleft Metodi Variazionali ed Equazioni Differenziali Non Lineari\textquotedblright }}
\author{Castelli Roberto, Terracini Susanna}
\date{}

\begin{document}
\maketitle
\vskip 20pt
\begin{abstract}
\noindent 

We study the possible regularization of collision solutions for one centre problems with a weak singularity. In the case of logarithmic singularities, we consider the method of regularization via smoothing of the potential. With this technique, we prove that the extended flow, where collision solutions are replaced with a transmission trajectory, is continuous, though not differentiable, with respect to the initial data.  
\bigskip


\end{abstract}

\section{Introduction}

In this paper we deal with dynamical systems associated with conservative central forces which are singular at the origin. A \emph{classical solution} does not interact with the singularity of the force, i.e., it is a path $u\in\mathcal{C}^2(T,\mathbb{R}^2\setminus\{0\})$ which fullfils  the initial value problem
\begin{equation}\label{sistema}
P\ :\ \left\{\begin{array}{ll}
&{\displaystyle\ddot{u}=\nabla V(|u|)}\\
& \\
& {\displaystyle(u(0),\dot{u}(0))=(q_0,p_0) \in (\mathbb{R}^2\setminus\{0\}\times\mathbb{R}^2)} \\
\end{array} \right.
\end{equation}
where $V(x)\in\mathcal{C}^2(\mathbb{R}^+,\mathbb{R})$ is the potential of the force and  $T$ denotes the maximal interval of existence.
As well-known, the two-body problem with an interaction potential $V$ can be reduced to a system of this form where  $u(t)$ denotes the position of one of the particle with respect to the center of mass. Accordingly, we shall term \emph{collision}  the configuration $u(t)=0$. As the force field diverges at $u=0$, collisions are among the main sources of non-completeness of the associated flow.  This work  studies the possible extensions of the flow through the collision that make it  continuous with respect to the initial conditions. We are concerned with weak singularities of the potential, namely logarithms. 

The regularization of total and partial collisions in the $N$-body problem is a very classical subject and, in the years, different strategies have been developed in order to extend motions beyond the singularity  \cite{LC,KS,Mo,Mg,E,SF,SM}. Very roughly, these classical methods rely upon suitable changes of space-time variables aimed at obtaining a smooth flow, possibly on an extended phase space;  to this aim, the first step is to determine the asymptotic behaviour of the collision solutions and then the phase space is extended either by means of a double covering, or with the attachment of a collision manifold.

In this paper we consider a further, non classical way of extending the flow, related to the technique of regularization via smoothing of the potential introduced by Bellettini, Fusco, Gronchi \cite{G}. This method consists in smoothing the singular potential and passing to the limit as the smoothing parameter $\varepsilon$ and the angular momentum tend to zero simultaneously but in an independent manner (indeed we know that the only collision motions have vanishing angular momentum). This involves an in-depth analysis about the ways the smoothing of the potential coupled with the perturbation of initial conditions lead to define a global solution of the singular problem. This technique, when successful, has the advantage of being extremely robust with respect to the application of existence techniques such as the direct method of the calculus of variation. Let us mention that variational methods have been widely exploited in the recent literature in order to obtain selected symmetric trajectories for $N$-body problems with Kepler potentials (\cite{FT}).

To begin with, we remove the singularity at $x=0$ and we denote with  $V_\varepsilon(x)$  the smoothed function defined as
$$
V_\varepsilon(x)=V(\sqrt{x^2+\varepsilon^2})\qquad \varepsilon>0
$$ 
Then we look at the regularized problem 
\begin{equation}\label{sistemareg}
P(\varepsilon)\ :\ \left\{\begin{array}{ll}
&{\displaystyle\ddot{u}=\nabla V_\varepsilon(|u|)} \\
& \\
& {\displaystyle(u(0),\dot{u}(0))=(q_0,p_0) \quad \in \mathbb{R}^2 \times \mathbb{R}^2\;.}
\end{array} \right.
\end{equation}
Unlike \eqref{sistema}, the differential equation \eqref{sistemareg} is no longer singular, so that the initial value problem admits a global solution in $\mathbb{C}^\infty((-\infty,+\infty);\mathbb{R}^2)$ for every choice of the datum $(q_0,p_0)$, provided $\nabla V(x)$ is sublinear at infinity\footnote{Without any additional assumption on the behaviour of the potential $V(x)$ far away from the origin, a solution of system \eqref{sistema} might have  singularities other than collisions: for instance solutions could blow up in finite time.}.
Since  we focus on the singularities due to collisions, we  fix a ball $B_0(\bar R)$ of radius $\bar R$ centered at the origin, where the collision is the only singularity that system \eqref{sistema} can develop and we denote with $S(V)\subset\mathbb{R}^2\times\mathbb{R}^2$  the set of initial conditions $(q_0,p_0)$ leading to collision for the  system $P$ with $|q_0|\leq \bar R$. For every $\bar \nu\in S(V)$ let $u_{\bar \nu}(t)\in\mathcal{C}^2(T,\mathbb{R}^2)$ be the  collision solution  where $T$ denote the maximal interval of existence such that $|u_{\bar \nu}(t)|\leq \bar R$. Denoting with $u_{\varepsilon,\nu}(t)$ the solution of \eqref{sistemareg} with initial data $\nu$, we investigate the existence of the asymptotic limit of the paths $u_{\varepsilon, \nu}(t)$ as $(\varepsilon,\nu)\rightarrow(0,\bar \nu)$, its relationship with the collision solution $u_{\bar \nu}(t)$ of the singular system $P$  and  the continuity of the limit trajectory with respect to initial data. The definition of regularization considered in \cite{G} is the following.
 
\begin{definizione}
Let $V(x)$ be a singular potential. We say that the problem  (\ref{sistema}) is weakly regularizable   via smoothing of the potential in $B_0(\bar R)$ if for every $\bar \nu\in S(V)$ there exist two sequences $(\varepsilon_k)_k$, $(\nu_k)_k$ tending to $0$ and $\bar \nu$ respectively, such that there exists
$$
\lim_{k\rightarrow\infty}u_{\varepsilon_k,\nu_k}=u_0 
$$ 
and the flow 
\begin{equation*}
\tilde u_\nu(t)=\left\{\begin{array}{ll}
u_\nu(t)\quad &\nu\not\in S(V)\\
u_0(t) &\nu\in S(V)
\end{array}\right.
\end{equation*}
is continuous with respect to $\nu$.
\end{definizione}
In addition we say that

\begin{definizione}\label{strongly}
The singular one centre  problem \eqref{sistema} is strongly regularizable via smoothing of the potential if there exists $\bar R$ such that for every $\bar \nu\in S(V)$ there exists
\begin{equation}\label{uniflim}
\lim_{(\varepsilon,\nu)\rightarrow(0,\bar \nu)}u_{\varepsilon,\nu}=u_0
\end{equation}
and the flow 
\begin{equation*}
\tilde u_\nu(t)=\left\{\begin{array}{ll}
u_\nu(t)\quad &\nu\not\in S(V)\\
u_0(t) &\nu\in S(V)
\end{array}\right.
\end{equation*}
is continuous with respect to $\nu$.
\end{definizione}

In both the definitions we mean that the limit of the regularizing paths $u_{\varepsilon,\nu}(t)$ and the continuity of the extended flow are held in the ball $B_0(\bar R)$.

In \cite{G} the authors prove that in the case of  homogeneous potential of degree $\alpha$, $V(x)=\frac{1}{|x|^\alpha}$, $\alpha>0$, the one-centre  problem is weakly regularizable via smoothing of the potential if and only if $\alpha$ is in the form 
\begin{equation}\label{alfa}
\alpha=2\left(1-\frac{1}{n} \right)
\end{equation}
where $n$ is a positive integer or $\alpha>2$. On the other hand  they show that the homogeneous problem is never strongly regularizable via smoothing of the potential. Indeed,  a  necessary condition in order to achieve the uniform limit \eqref{uniflim} is that the apsidal angle $\Delta\theta_l(u)$ of a solution of the system \eqref{sistema} has to converge to $\frac{\pi}{2}$ as the angular momentum $l$ tends to zero (see the definition of apsidal angle in the next section). This condition is never satisfied by  $\alpha$-homogeneous potentials $\alpha>0$, since $\Delta\theta_l(u)\rightarrow \frac{\pi}{2-\alpha}>\frac{\pi}{2}$ as $l\rightarrow 0$ \cite{G,TT}. This explains why the one centre problem with  homogeneous potential  can not be regularized according to definition \ref{strongly}. Conversely, when the logarithmic potential is considered, it can be proved (\cite{TT}) that the limiting apsidal angle do indeed converge to $\pi/2$. Then there is no obstruction and we could expect that the limit \eqref{uniflim} is attained. This fact suggests to extend the motion after a collision by reflecting it about the origin.
We will show that, in this way, not only for the logarithmic potential, but for a larger class $\mathcal{V}^*$ of potentials, the problem is regularizable according to definition \ref{strongly}.  

The sets of potential functions we will consider in this paper are the following.
\begin{definizione}[The function set $\mathcal{V}$ ]\label{defV}
We define $\mathcal{V}$ the set of functions $V(x)\in\mathbb C^\infty(\mathbb R^+, \mathbb R)$ with the properties:

\begin{minipage}{9cm}
\begin{itemize}
\item[\textbf{i.}] $\displaystyle \lim_{x\rightarrow 0^+}V(x)=+\infty$
\end{itemize}
there exists $S>0$ such that for every $x\in(0,S)$
\begin{itemize}
\item[\textbf{ii.}] $\displaystyle V'(x)<0$ ,  $V''(x)>0$ 
\item[\textbf{iii.}] the function $\displaystyle\frac{V'(x)}{V''(x)}$ is decreasing with respect to x
\end{itemize}
and
\begin{itemize} 
\item[\textbf{iv.}] $\displaystyle\frac{d}{dx}\frac{V'(x)}{V''(x)}(0)<-\frac{1}{2}$
\end{itemize}
\end{minipage}
\begin{minipage}{4cm}
\includegraphics[width=4cm]{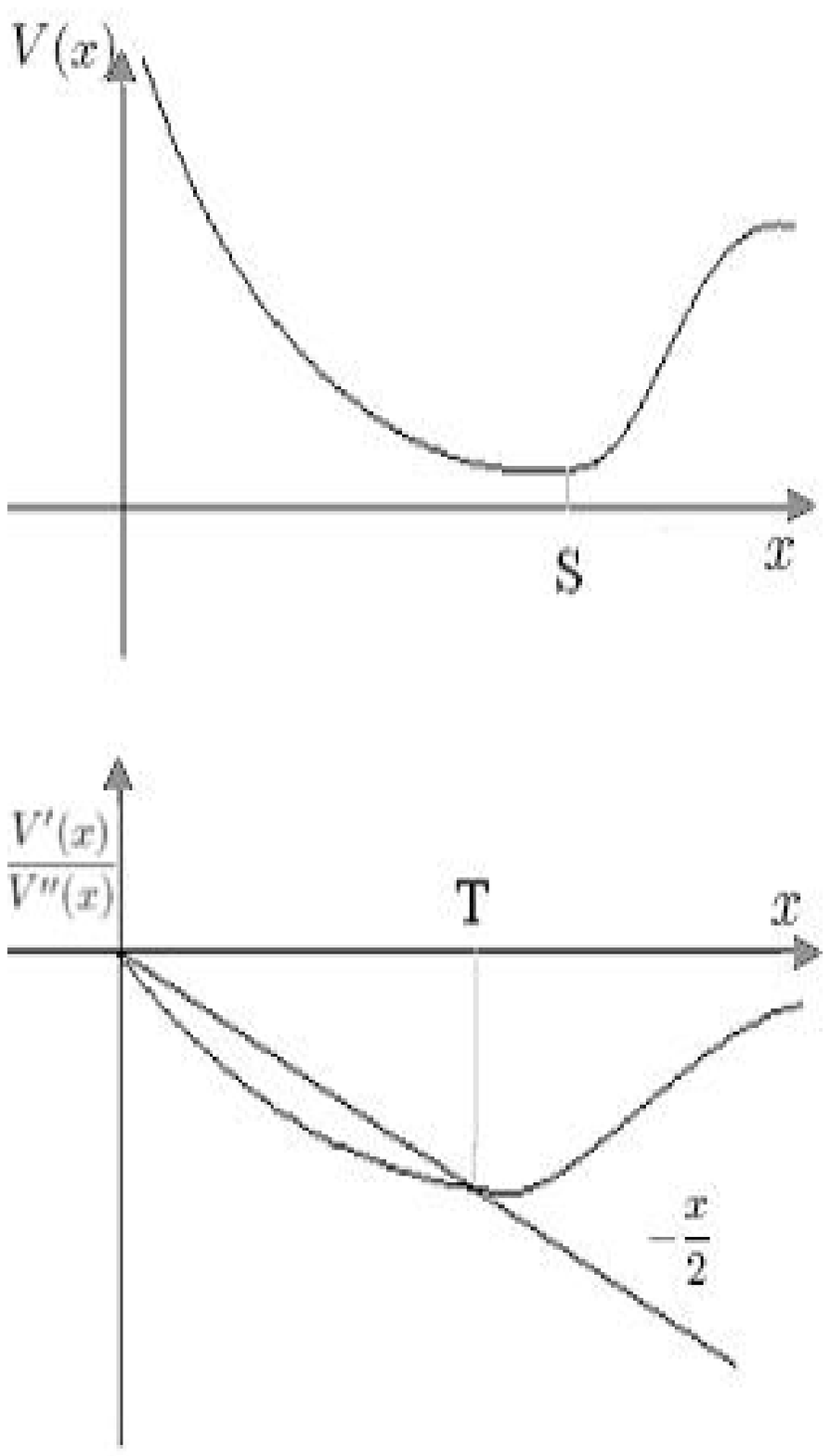}
\end{minipage}  
\end{definizione}

The properties \textbf{iii,iv} guarantee the existence of a $ T>0$  such that
\begin{equation}\label{definizT}
\displaystyle{
\frac{V'(x)}{V''(x)}\leq-\frac{x}{2} \quad \rm{in}\quad (0, T)}
\end{equation}
Let 
\begin{equation}\label{rvalidi}
\bar R:=\min\{ T,S\}.
\end{equation}

\begin{definizione}[The set $\mathcal{V}^*$]\label{defVstar}
Denote with $\mathcal{V}^*$ the set of functions $V(x)\in \mathcal{V}$ with the further property
\begin{itemize}
\item[\textbf {v.}]  $\displaystyle{\lim_{\lambda\rightarrow 0}\frac{V(\lambda x)}{V(\lambda)}=1}$ for every $x\geq 1 $ uniformly in every compact $K=[1,M]$. 
\end{itemize}
\end{definizione}

The set $\mathcal{V}$ includes potentials having homogeneous singularities and weaker. For instance the logarithmic potential, $V(x)=-\log(x)$,  as well as the homogeneous potentials, $V(x)=|x|^{-\alpha}$, provided $\alpha\in(0,1)$, belong to $\mathcal{V}$. On the other hand condition {\it\textbf{v.}} can be considered as  a logarithmic type property or a zero-homogeneity property: indeed it is never satisfied by homogeneous potentials, while the logarithmic potential is a prototype of all the functions satisfying condition {\it\textbf{v.}}.

Our main goal is the following:

\begin{teorema}\label{Mainteo}
For every $V(x)\in\mathcal{V}^*$ the one centre problem is regularizable according to definition \ref{strongly} where $\bar R$  is given in \eqref{rvalidi}.
\end{teorema}

In the particular case of logarithmic potential, $V(x)=-\log(x)$, one has $\bar R=+\infty$, therefore 

\begin{corollario}
The logarithmic one central problem is globally regularizable via smoothing of the potential according to definition \ref{strongly}.
\end{corollario}

The paper is organized as follows. In section \ref{secpreliminari} we follow the classical method for dealing with central problem based on first integrals and we derive the set $S(V)$ of initial conditions leading to the collisions for the unperturbed system. Next, in section \ref{setting}, given any collision solution $u(t)$, we set the initial data $\bar \nu\in S(V)$ and we define the family of paths $u_{\varepsilon,\nu}(t)$. Section \ref{dimostrazione} contains the proof of the main theorem and the  analysis of the regularity of the extended flow. The main part of the proof consists in proving the existence of the limit of the path $u_{\varepsilon,\nu}(t)$ as $(\varepsilon,\nu)\rightarrow(0,\bar \nu)$, especially for what that concerns  the angular part, (Theorem \ref{teoregolariz}). This is the most delicate step, for the it involves the uniformity of the limit as $(\varepsilon,\nu)\rightarrow(0,\bar \nu)$, and it allows to conclude the strong regularizability of the  problem. 

It results that the natural extension of the collision solution is the  \textit{transmission solution}, see definition \ref{def-trasm}, obtained by reflecting the motion through the collision.  The regularity of the extended flow is carried on in section \ref{secreg}: in theorems \ref{teopoincare} and \ref{teocontiuitapoincare} the continuity  of the Poincar\'e  map and Poincar\'e  section with respect initial data is achieved.

Furthermore, in order to have a complete picture of the problem, in section \ref{secvar} we join a variational approach and we analyse the variational properties of the collision paths.

\section{Preliminaries}\label{secpreliminari}
For any choice of the potential function $V(x)\in\mathcal{C}^2(\mathbb{R}^+,\mathbb{R})$ the one centre problem \eqref{sistema} is a hamiltonian system and admits two integrals of motion: the energy $E$ and the angular momentum $l$:
$$
E=\frac{1}{2}|\dot u|^2-V(|u|)\qquad l=\dot u\wedge u
$$
Since the conservation of the angular momentum implies the motion is planar, in the following $l$ is used to denote the modulo of the angular momentum, rather than the vector.
The radial symmetry of the equation of motion \eqref{sistema}  suggests to introduce the polar coordinates in the plane $(r,\theta)$. In this setting  the quantities $E$ and $l$ are expressed in the form
\begin{equation}\label{formulaenergia}
E=\frac{1}{2}\dot r ^2+\frac{1}{2}\frac{l^2}{r^2}-V(r)\qquad 
l=r^2\dot\theta
\end{equation}
We define
\begin{equation}\label{funzionef}
f(r)=2r^2(E+V(r))
\end{equation}
 then the  relation \eqref{formulaenergia} reads as $
l^2=f(r)-(r\dot r)^2$
\begin{figure}[t]
\subfigure[Apsidal values]{
\includegraphics[scale=0.48]{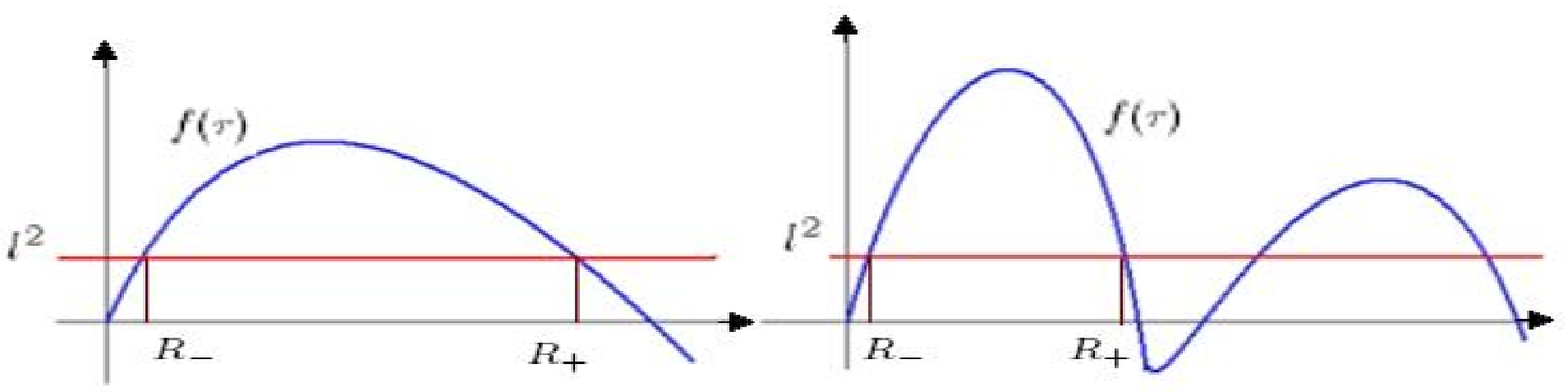}\label{figapsidalval}}
\subfigure[Apsidal angle]{
\includegraphics[scale=0.22]{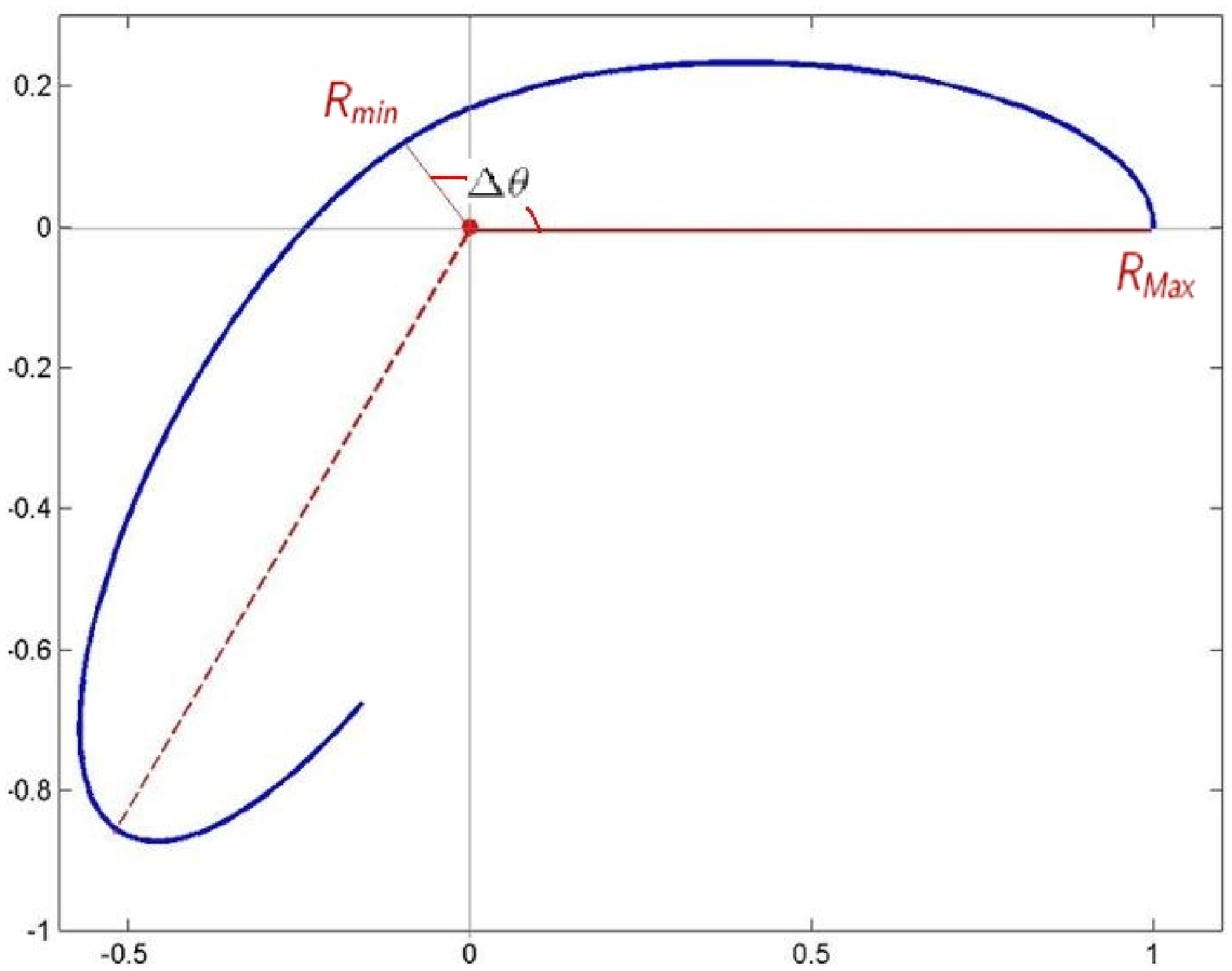}\label{figapsidalan}}
\caption{}
\end{figure}
and  shows that a solution of system (\ref{sistema}) of energy $E$ and angular momentum $l$ exists only for those values of the radial coordinate $r\geq 0$ such that $f(r)\geq l^2$.

For fixed values of $E$ and $l$, we denote with $R_+$, if it exists, the minimum positive value of $r$ such that $f(r)=l^2$ and $f'(r)<0$, $R_+=+\infty$ otherwise,  and we denote with $R_-$, if it exists, the minimum positive value of $r<R_+$  such that $f(r)=l^2$ and $f'(r)>0$, $R_-=0$ otherwise, see figure \ref{figapsidalval}.  By definition, it follows that $R_-=0$ for collision solutions and $R_+=+\infty$ for unbounded orbits. In any case  we term $R_+$ and $R_-$  the maximal and minimal value of the angular coordinate and we refer to them as the\textit{ apsidal values} of the orbit. Moreover, following the terminology adopted in celestial mechanics, we sometimes refer to $R_+$ and $R_-$ respectively  as the apocentre and the pericentre of the orbit. As it is well known, in case of  non collision and bounded trajectories, the radial coordinate $r(t)$ oscillates   periodically between its extremal values $R_+$ and $R_-$ while the angular coordinate $\theta(t)$ covers an angle equal to 
$$
\displaystyle\Delta\theta_l(u)=\displaystyle{\int_{R_-}^{R_+}
\frac{1}{r\sqrt{\frac{2r^2}{l^2}(E+V(r))-1}}dr}
$$   
between  each singular oscillation of $r(t)$.
We term the angle $\Delta\theta_l(u)$ the \textit{apsidal angle}, see figure \ref{figapsidalan}.

 The knowledge of the apsidal values $R_+$ and $R_-$ and the value of the apsidal angle is sufficient to determine the behaviour of the solution since the whole trajectory is obtained repeating periodically the part of path between a point where $r(t)$ is maximum and the following point where  $r(t)$ is minimal. 

The definition of the apsidal angle extends in a natural way for unbounded and collision solutions:
in the first  case the orbits is composed by a single oscillation of the radial coordinate from infinity to $R_-$ and back to infinity and the apsidal angle represents the angle covered by the particle coming from infinity to the pericentre and it is obtained replacing $R_+=+\infty$ in the previous relation. On the other hand, if a collision occurs, the apsidal angle denotes the increment of the angular coordinate between the apocentre and the collision point.

In order to characterise the set $S(V)$ of initial data leading to a collision  we give the following definition.  
\begin{definizione}\label{defstrong}
We say that a potential function $V(x)$, singular in the origin, is  of {weak type} if 
$$
\lim_{x\rightarrow 0^+}x^2V(x)=0
$$
Otherwise we say that $V(x)$ is a {strong type} potential.
\end{definizione}
A similar classification of singular potentials can be found in a work of Gordon \cite{Go} where a potential $V(x)$ is said to satisfy a \textit{strong force condition} at a point $x_0$ if $V(x)$ tends to infinity as $x$ tends to $x_0$ and also there exists a function $U(x)$ with infinitely deep wells at $x_0$, such that
$$
V(x)\geq|\nabla U(x)|^2
$$
in a neighbourhood of $x_0$. We  notice that, among the homogeneous potentials, the set of potentials with the property to be of strong type and the ones  satisfying  the Gordon's  strong force condition coincide. Clearly, an $\alpha$-homogeneous potential is of weak type if and only if $\alpha\in(0,2)$ and in these cases a  collision occurs  only in zero angular momentum orbits \cite{Mg}, while if $\alpha\geq 2 $ a collision solution exists also for non-zero values of the angular momentum \cite{G,Mg}. 
The next proposition extends this result.

\begin{proposizione}\label{teoweaktype}
If $V(x)\in\mathcal{C}^2(\mathbb{R}^+,\mathbb{R})$ is a potential of weak type, a solution $u(t)$ of  the dynamical system  $\ddot u=~\nabla V(|u|)$ ends into a collision if and only if the angular momentum is zero.
\end{proposizione}
\dimostr 
Denote with $E$ and $l$ the energy and the angular momentum of the solution $u(t)$ and let $f(r)$ as in \eqref{funzionef}. As mentioned before, a solution 
exists only for the values of radial coordinate $r$ satisfying $
l^2\leq f(r)$.
Suppose $l=0$: since $V(x)$ tends to infinity as $x$ goes to zero, for every value of $E$ there exists  a neighbourhood of the origin
 where $E+V(r)>0$ then the solution presents a collision.

Conversely, since $V(x)$ is a weak type potential, it follows that $f(r)\rightarrow 0$ as $x\rightarrow 0^+$ thus for every value of $l\not=0$  
there exist a neighbourhood of the origin where $l^2>f(r)$. Hence the collision can not be attained on solutions with non zero angular momentum.
\qed

\begin{proposizione}\label{teodicollisione}
Every $V(x)\in\mathcal{V}$ is a weak type potential.
\end{proposizione}
\dimostr
From relation \eqref{definizT}, by integration, it  follows the estimate
$$
-V'(\xi)\leq\frac{C_1}{\xi^{2}}\ ,\qquad C_1> 0
$$
for every $\xi\in(0,T)$. Therefore, again by integration, we infer
\begin{equation}\label{stimaV}
V(x)\leq\frac{C_1}{x}+C_2
\end{equation}
and we conclude
$$
\lim_{x\rightarrow 0}x^2V(x)=0\ .
$$
\qed

Propositions \ref{teoweaktype}, \ref{teodicollisione} show that, for any choice of $V(x)\in\mathcal{V}$, a solution of the system \eqref{sistema} ends into a collision  if and only if the angular momentum is zero. Therefore the set $S(V)$ of initial conditions $\bar \nu=(\bar q_0,\bar p_0)$ that make the solutions to be singular consists in those $(\bar q_0,\bar p_0)$ satisfying $l=|\bar q_0\wedge \bar p_0|=0$.

\section{Setting}\label{setting}
For every fixed  $V(x)\in\mathcal{V}$ let $\bar R$ be the quantity defined in \eqref{rvalidi} and let  $B_0(\bar R)$ be used to  denote  the ball of radius $\bar R$ around the origin. The properties of the potential class $\mathcal V$ assure that  the collision is the only source of singularity for the dynamical system inside $B_0(\bar R)$. As discussed in the introduction, given a collision solution $\bar u(t)$ for the one central problem \eqref{sistema}, our intent is to define an extension in $B_0(\bar R)$ for the solution $\bar u(t)$ beyond the collision.

To this we first have  to set the initial conditions $\bar \nu=(q_0,p_0)$  for the singular path $\bar u(t)$:
we denote with $P$ the first positive solution  of equation $f(r)=0$, $P:=+\infty$ if such value does not exists. We remind that the angular momentum $l$ is zero for collision solution, hence the orbits drawn by $\bar u(t)$ is a straight line joining some point in the plane with the origin.   

An alternative occurs:

\begin{figure}[h]
 \begin{minipage}[b]{6.5cm}
   \centering
   \includegraphics[width=6cm]{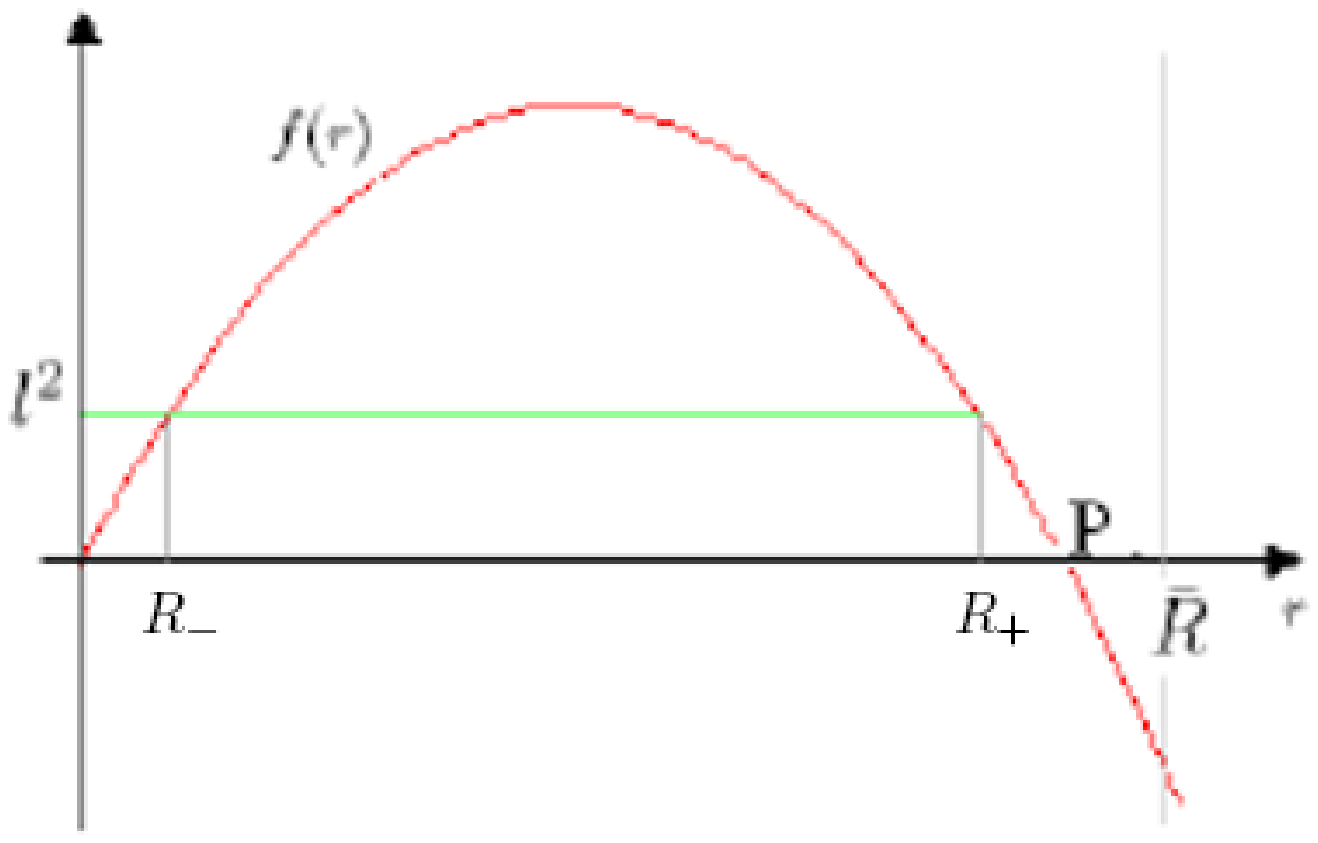}
 \label{figura1a}
 \end{minipage}
 \ \hspace{2mm} \hspace{3mm} \
 \begin{minipage}[b]{6.5cm}
  \centering
   \includegraphics[width=5cm]{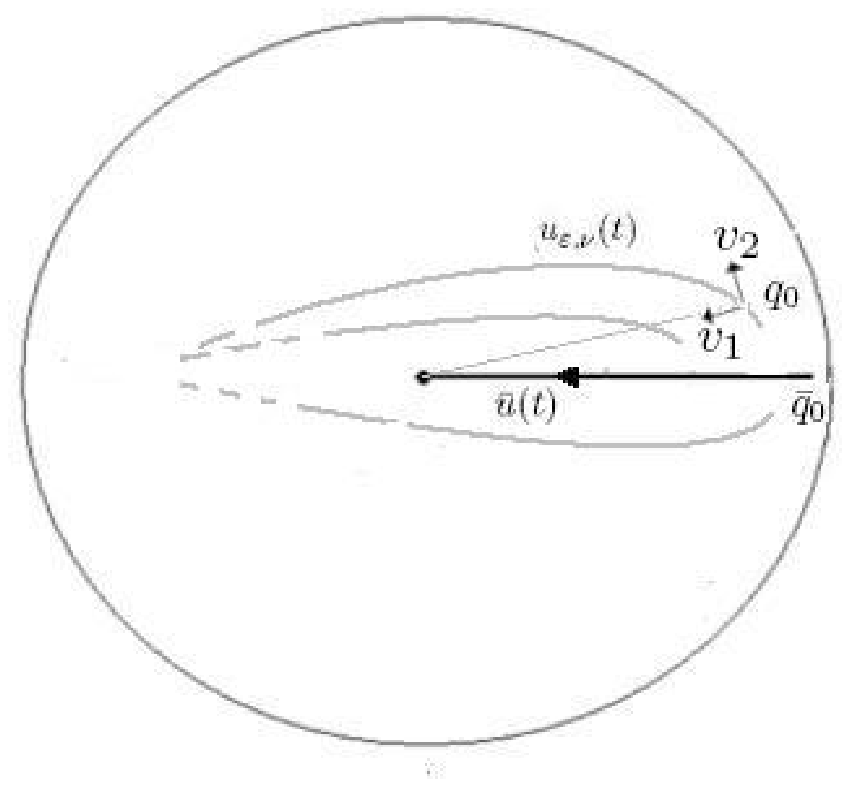}
\label{figura1b}

 \end{minipage}
\caption{Case 1 }
\label{figura2}
\end{figure}

\noindent \textbf{Case 1} $P<\bar R$.  

The collision solution $\bar u(t)$ is bounded in a ball centered in the origin of radius $P<\bar R$ and  $P$ represents the maximal value of the radial coordinate ( Figure \ref{figura2}). Without loss of generality, we can set the initial condition $\bar\nu$  of the collision solution $\bar u$ as 
\begin{equation}\label{iniz1coll}
\bar \nu=(\bar q,0), \qquad|\bar q|=P
\end{equation}

\vskip 5pt

\noindent\textbf{Case2} $P\geq \bar R$

\begin{figure}[h]
 \begin{minipage}[]{6.5cm}
   \centering
   \includegraphics[width=6cm]{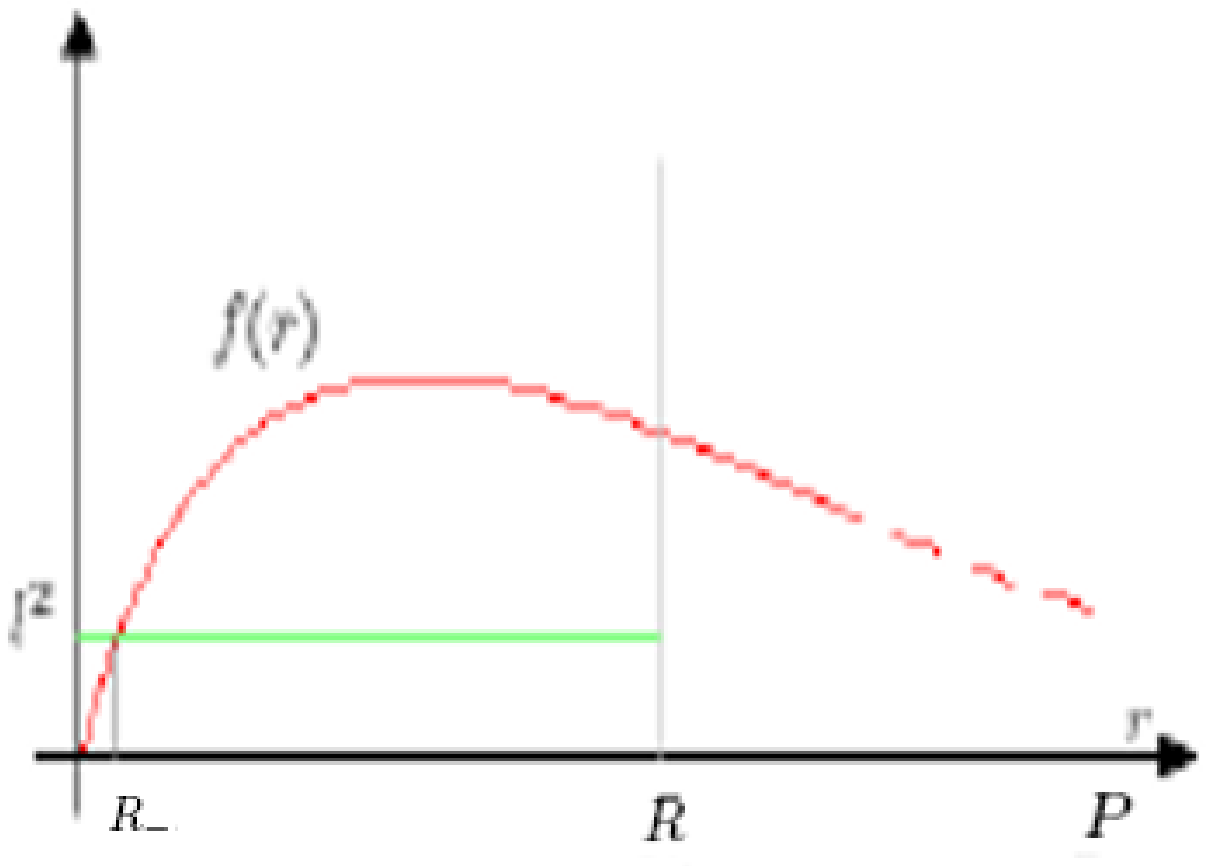}
\label{figura2a}
 \end{minipage}
 \ \hspace{2mm} \hspace{3mm} \
 \begin{minipage}[]{6.5cm}
  \centering
   \includegraphics[width=5cm]{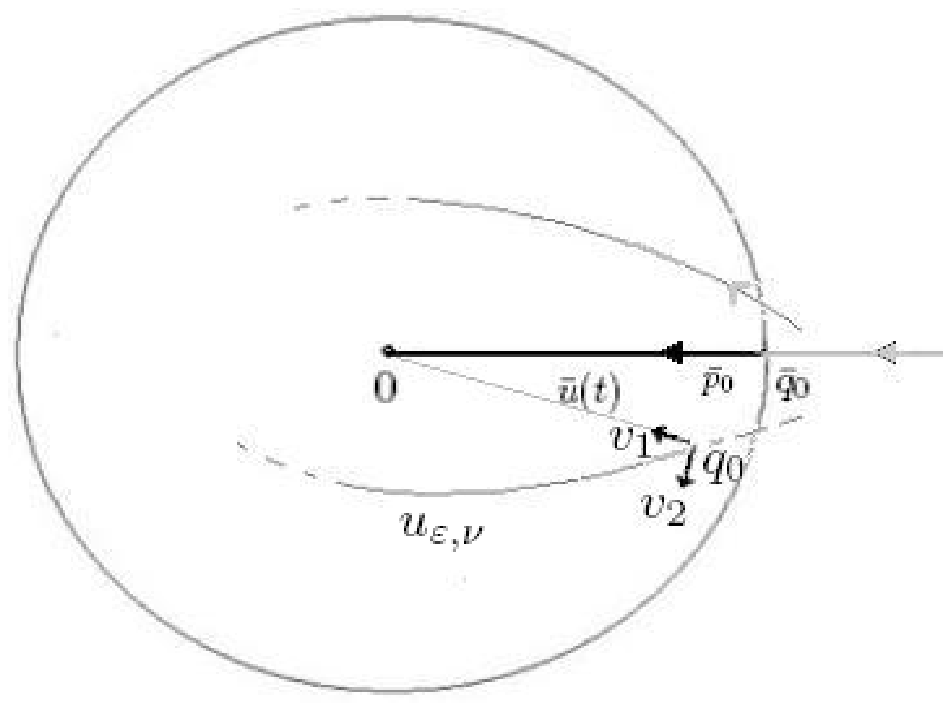}
\label{figura2b}
 \end{minipage}
\caption{Case 2}
\label{figura3}
\end{figure}

In this case the collision solution is not bounded in $B_0(\bar R)$ and it could  also be unbounded. We focus our analysis only on the portion of path bounded by $\bar R$ hence we select as initial condition for $\bar u(t)$  the couple 
\begin{equation}\label{iniz2coll}
\bar \nu=(\bar q,\bar p),\qquad|\bar q|=\bar R,\qquad |\bar p|^2=2(E+V(\bar R))
\end{equation}
where the initial velocity $\bar p$ is directed toward the center of attraction ( Figure \ref{figura3}).

In both the cases, for  $\varepsilon>0$ let $u_{\varepsilon,\nu}(t)$ be a  solution of system \eqref{sistemareg} leading from an initial data  $\nu=(q_0,p_0)$.
We refer to $u_{\varepsilon,\nu}(t)$ as \textit{regularizing paths} in order to underline the purpose to define the extension for the singular solution $\bar u(t)$ as a limit of $u_{\varepsilon,\nu}(t)$ as $({\varepsilon,\nu})\rightarrow(0,\bar \nu)$. The smoothing of the potential does not affect the hamiltonian structure of the system, therefore the angular momentum $l$ and the energy $E_{\varepsilon,\nu}=\frac{1}{2}|\dot u_{\varepsilon,\nu}|^2-V_\varepsilon(|u_{\varepsilon,\nu}|)$ are conserved along the solution $u_{\varepsilon,\nu}(t)$.
Decomposing the initial velocity $p_0$ in terms of the parallel and orthogonal component with respect to $q_0$,
\begin{equation}
p_0=v_1+v_2\qquad v_1\parallel q_0,\quad  v_2\perp q_0
\end{equation} one has
$|v_2|^2=\frac{l^2}{|q_0|^2}$ and $E_{\varepsilon,\nu}=\frac{1}{2}|v_1|^2+
\frac{1}{2}\frac{l^2}{|q_0|^2}-V_\varepsilon(|q_0|)$ then the condition $\nu\rightarrow \bar \nu$ is  equivalent to 
$$
q_0\rightarrow\bar q,\quad l\rightarrow 0,\quad  v_1\rightarrow \bar p 
$$
that, coupled with the condition $\varepsilon\rightarrow 0$, implies $E_{\varepsilon,\nu}\rightarrow E$. It turns out that also the apsidal values of the regularizing paths have to converge to the corresponding ones of $\bar u(t)$, indeed 
as $(\varepsilon,\nu)\rightarrow(0,\bar \nu)$,  the pericentre $R_-$ of the solution $u_{\varepsilon,\nu}(t)$ tends to zero while the apocentre $R_+$ is bounded by $\bar R$ and tends to $P$ in \textbf{case 1}, while $R_+>\bar R$ and possibly $R_+=+\infty$ in \textbf{case 2}.

In the following sections we will  
deal with the existence and the property of the limit for the paths $u_{\varepsilon,\nu}(t)$ as ${(\varepsilon,\nu)}\rightarrow(0,\bar \nu)$. As discussed in \cite{G}, the  behaviour of the angular coordinate of the regularizing paths plays a fundamental role for the existence of the  limit of $u_{\varepsilon,\nu}(t)$ uniformly as  ${(\varepsilon,\nu)}\rightarrow(0,\bar \nu)$ rather then for subsequences 
${(\varepsilon_k,\nu_k)}\rightarrow(0,\bar \nu)$.
For this reason and since we focus our analysis only inside the ball $B_0(\bar R)$, we extend the definition of apsidal angle for the solution $u_{{\varepsilon,\nu}}(t)$ as follows:  we denote with $\Delta\theta(u_{\varepsilon,\nu})$ the apsidal angle of  the path $u_{\varepsilon,\nu}(t)$ in case the trajectory is  bounded by $\bar R$, \textbf{case1}
\begin{equation}\label{dteta1}
\Delta\theta(u_{\varepsilon,\nu})=\int_{R_-}^{R_+}\frac{l}{r^2\dot r}dr
\end{equation}

otherwise, in \textbf{case 2}, we denote with  $\Delta\theta(u_{\varepsilon,\nu})$ the angle covered by the path $u_{\varepsilon,\nu}(t)$ between the point where $u_{\varepsilon,\nu}(t)$ enter in the ball $B_0(\bar R)$ and the point of minimal distance from the origin 
 \begin{equation}\label{dteta2}
\Delta\theta(u_{\varepsilon,\nu})=\int_{R_-}^{\bar R}\frac{l}{r^2\dot r}dr\ .
\end{equation} 
 
\section{Proof of Main Theorem and property of the extended flow}\label{dimostrazione}
The proof of  theorem \ref{Mainteo} is composed by two parts: first, in section \ref{secproof} we prove the existence of the limit of the trajectories  $u_{\varepsilon,\nu}(t)$ as $(\varepsilon,\nu)\rightarrow(0,\bar \nu)$  and we define the extension of the singular solution, then in section \ref{secreg} we study the regularity of the extended flow.  

A necessary condition for the existence of the limit \eqref{uniflim} is the existence of the limit of the angular part of the regularized solutions. Theorem \ref{teoregolariz} concerns the asymptotic of $\Delta\theta(u_{\varepsilon,\nu})$ as $({\varepsilon,\nu})\rightarrow(0,\bar \nu)$: to  this aim we first prove in lemma \ref{limitatezzak1k2} the $\mathbb{L}^1$ boundness of the integrand in \eqref{dteta1} and \eqref{dteta2} then we apply the dominated convergence theorem and pass to the limit under the integral sign. 
The boundness of the integrand is a consequence of a technical estimate stated in the proposition  \ref{limitatezzaK} and it is attained for every potential $V(x)\in\mathcal V$, while the existence of the limit is a consequence of proposition \ref{allimite} based on assumption \textbf {v.}.
The result we obtain suggest to define the extension $u_0(t)$ of the collision solution $\bar u(t)$ beyond the singularity as a \textit{transmission solution}.

In order to gain the regularity of the extension, we analyze, in theorems \ref{teopoincare} and \ref{teocontiuitapoincare}, the continuity  of the  Poincar\'e map and Poincar\'e section of the extended flow in the phase space. 

\subsection{The existence of the limit of $u_{\varepsilon,l}(t)$}\label{secproof}
\label{teoremi di limitatezza}
\begin{proposizione}\label{limitatezzaK}
Let $V(x)\in\mathcal{V}$  and $\bar R$ as in \eqref{rvalidi}. Then for every $\bar r<\bar R$ there exists an $\bar\varepsilon$  such that $\forall \varepsilon<\bar \varepsilon$ and for every $0<y<\bar r$ it holds
\begin{equation}\label{limitato}
F_{(y,\bar r)}(\varepsilon,x):=\frac{(\bar r+x)}{(\frac{x}{y}-1)}\left[\left(\frac{x}{y}\right)^2\left(
\frac{V_\varepsilon(x)-V_\varepsilon(\bar r)}{V_\varepsilon(y)-V_\varepsilon(\bar r)}\right)
\frac{(\bar r^2-y^2)}{(\bar r^2-x^2)}-1\right]\geq \bar r
\end{equation}
for every $y\leq x\leq \bar r$.
\end{proposizione}
The proof of the proposition \ref{limitatezzaK} is split in two parts: we first show that there exists $\bar \varepsilon>0$ such that for every $\varepsilon<\bar \varepsilon$ the function $F_{(y,\bar r)}(\varepsilon,x)>F_{(y,\bar r)}(0,x)$ for every $0<y\leq x\leq \bar r<S$ , then we show that $F_{(y,\bar r)}(0,x)\geq\bar r $ for every  $0<y\leq x\leq \bar r<\bar R$. 
\begin{lemma}\label{limiteepsilon}
Let $V(x)$ be a function satisfying  the properties \textbf{i.-iii.}. Then for every choice of $0<~y<~\bar r<~S$ there exists $\bar\varepsilon >0$ such that $\forall \varepsilon<\bar\varepsilon$ and $\forall x\in(y,\bar r)$ it holds
\begin{equation}\label{limiteK}
Q_{(y,\bar r)}(\varepsilon,x):=\frac{V_\varepsilon(x)-V_\varepsilon(\bar r)}{V_\varepsilon(y)-V_\varepsilon(\bar r)}\geq \frac{V(x)-V(\bar r)}{V(y)-V(\bar r)}
\end{equation}
\end{lemma}
\dimostr
By means of straightforward calculations and reminding the definition of smoothed potential $V_\varepsilon(\cdot)=V(\sqrt{(\cdot)^2+\varepsilon^2})$, the relation (\ref{limiteK}) is equivalent to 

\begin{equation*}
\frac{V(\sqrt{y^2+\varepsilon^2})-V(\sqrt{x^2+\varepsilon^2})}{V(\sqrt{x^2+\varepsilon^2})-V(\sqrt{\bar r^2+\varepsilon^2})}  \leq\frac{V(y)-V(x)}{V(x)-V(\bar r)}
\end{equation*}
For every $s>0$ we define the function $U(s)=V(\sqrt s)$. Obviously the function $U(s)$ inherit property \textbf{i.} and property \textbf{ii.} for every $s\in(0,S^2)$, while the relation
$$
2\sqrt{s}\frac{U''(s)}{U'(s)}=\frac{V''(\sqrt{s})}{V'(\sqrt{s})}-\frac{1}{2\sqrt{s}}
$$
and property \textbf{iii.} implies that the function $\displaystyle{\sqrt{s}\frac{U''(s)}{U'(s)}}$ is increasing for every $s\in(0,S^2)$.
The proof of the  lemma follows proving that for every choice of $0<~y<~\bar r<~S^2$ there exists $\bar\varepsilon >0$ such that $\forall \varepsilon<\bar\varepsilon$ and $\forall x\in(y,\bar r)$ it holds
\begin{equation}
\frac{U(y+\varepsilon)-U(x+\varepsilon)}{U(x+\varepsilon)-U(\bar r+\varepsilon)}\leq\frac{U(y)-U(x)}{U(x)-U(\bar r)}.
\end{equation}
Let   $0<y<\bar r<S^2$ be fixed and let the function $g(\varepsilon,x)$ be defined
as 
$
g(\varepsilon,x):=\frac{U(y+\varepsilon)-U(x+\varepsilon)}{U(x+\varepsilon)-U(\bar r+\varepsilon)}$, $ x\in(y,\bar r)$.
Since $g(0)=\frac{U(y)-U(x)}{U(x)-U(\bar r)}$  it's enough to show that 
\begin{equation}\label{dg0}
\frac{dg}{d\varepsilon}(0,x)<0\quad \forall\ 0<y<x<\bar r<S^2\ .
\end{equation} 
The sign of the derivative
is given by the sign of the function 
$$
G_{(y,\bar r)}(x):=U'(y)\Big(U(x)-U(\bar r)\Big)+U(y)\Big(U'(\bar r)-U'(x)\Big)+U'(x)U(\bar r)-U'(\bar r)U(x)
$$
 We observe  that $
G_{(y,\bar r)}(y)=G_{(y,\bar r)}(\bar r)=0
$
and
\begin{align}\label{Dg}
\frac{dG_{(y,\bar r)}}{dx}(x)=U'(x)\Big(U'(y)-U'(\bar r)\Big)+U''(x)\Big(U(\bar r)-U(y)\Big)
\end{align}
Since $G_{(y,\bar r)}(x)$ is continuous,  there exists at least one point $\bar x\in(y,\bar r)$ where $\frac{dG_{(y,\bar r)}}{dx}(\bar x)~=~0$; the proof of the lemma follows once we prove the inequalities
\begin{equation}\label{disug}
\begin{array}{ll}
 {\displaystyle \frac{dG_{(y,\bar r)}}{dx}(x)<0}\qquad {\rm for }\quad x\in (y,\bar x)\\
 {\displaystyle \frac{dG_{(y,\bar r)}}{dx}(x)>0}\qquad {\rm for }\quad x\in (\bar x,\bar r)\\
\end{array}
\end{equation}
Indeed, if \eqref{disug} hold, the function $G_{(y,\bar r)}(x)$ and, consequently the derivative in \eqref{dg0} are negative for every $x\in(y,\bar r)$. To obtain the relations \eqref{disug} we multiply both the sides in \eqref{Dg} for $\sqrt{x}$ and divide them for $U'(x)(U(\bar r)-U(y))$. We obtain 
$$
N(x):=\frac{\sqrt x G'_{y,\bar r}(x)}{U'(x)(U(\bar r)-U(y))}=\sqrt x\frac{U''(x)}{U'(x)}+\sqrt x\frac{(U'(y)-U'(\bar r))}{(U(\bar r)-U(y))}
$$
By definition of $\bar x$ it follows $N(\bar x)=0$. Moreover, since 
$\displaystyle{\sqrt{s}\frac{U''(s)}{U'(s)}}$ is increasing for every $s\in(0,S^2)$ and the factor $\displaystyle{\frac{(U'(y)-U'(\bar r))}{(U(\bar r)-U(y))}}$ is positive, we infer that $N(x)$ is increasing in $x$.
Therefore $N(x)<0$ for $x<\bar x$ and $N(x)>0$ otherwise and, since $U'(x)(U(\bar r)-U(y))>0$, inequalities \eqref{disug} hold. 
\qed
\vskip 5pt
\noindent {\it Proof of Proposition \ref{limitatezzaK}.}

\noindent We fix $0<y<\bar r<\bar R$. For  lemma \ref{limiteepsilon} there exists $\bar \varepsilon>0$ such that  $F_{(y,\bar r))}(\varepsilon,x)~\geq~ F_{(y,\bar r))}(0,x)$ for every $\varepsilon\in[0,\bar \varepsilon]$ and for every $x\in[y,\bar r]$.
Hence  it's sufficient to prove  
\begin{equation*}
F_{(\bar r,y)}(x):=\frac{(\bar r+x)}{(\frac{x}{y}-1)}\left[\left(\frac{x}{y}\right)^2\left(
\frac{V(x)-V(\bar r)}{V(y)-V(\bar r)}\right)
\frac{(\bar r^2-y^2)}{(\bar r^2-x^2)}-1\right]\geq\bar r
\end{equation*}
Suppose for a moment that relation
\begin{equation}\label{condizione}
\frac{V(x)-V(\bar r)}{V(y)-V(\bar r)}\geq \left(\frac{\bar r -x}{\bar r -y}\right)\frac{y}{x}
\end{equation}
holds for every $x\in(y,\bar r)$.
Replacing into $F_{(\bar r,y)}(x)$ we obtain
\begin{align*}
F_{(\bar r,y)}(x)\geq\frac{(\bar r+x)}{(\frac{x}{y}-1)}\left[\left(\frac{x}{y}\right)\frac{(\bar r+y)}{(\bar r+x)}-1\right]\geq \bar r\ .
\end{align*}

In order to prove  relation (\ref{condizione}) we rewrite it as
\begin{equation}\label{relazione}
\frac{(V(x)-V(\bar r))}{(\bar r-x)}x\geq\frac{(V(y)-V(\bar r))}{(\bar r-y)}y \quad \forall x\in(y,\bar r)
\end{equation}
For  $x=y$ the inequality holds; moreover, denoting with $N(x)$    the numerator of the derivative
\begin{gather}
\frac{d}{dx}\left(\frac{(V(x)-V(\bar r))}{(\bar r-x)}x\right)=
\frac{x(\bar r-x)V'(x)
+\bar r (V(x)-V(\bar r))}
{{(\bar r-x)^2}}\label{derivata1}
\end{gather}
one has $N(\bar r)=0$ and $ 
\frac{dN}{dx}(x)=(\bar r -x)\Big(2V'(x)+xV''(x)\Big).
$
For  \eqref{definizT}, for every $x\in(0,\bar R)$,  $  \frac{dN}{dx}(x)\leq 0$, hence the derivative in (\ref{derivata1}) is positive and relation (\ref{relazione}) holds for every $x\in(y,\bar r)$.
\qed

\begin{proposizione}\label{allimite}
Let $V(x)\in\mathcal{V}^*$, then for every $\rho>1$
$$
\lim_{(\delta,\varepsilon)\rightarrow(0,0)}
\frac{V_\varepsilon(\rho \delta)}{V_\varepsilon(\delta)}=1 \qquad \delta,\varepsilon>0
$$
\end{proposizione}
\dimostr
We rewrite the above limit in the form
$$
\lim_{(\delta,\varepsilon)\rightarrow(0,0)}\frac{V
\Big(\sqrt{\delta^2+\varepsilon^2}\sqrt{\frac{\rho^2\delta^2+\varepsilon^2}
{\delta^2+\varepsilon^2}}
\Big)}{V(\sqrt{\delta^2+\varepsilon^2})}.
$$
Since $\rho>1$, for every choice of positive values of $\varepsilon$ and $\delta$,  it holds
$$
1\leq\frac{\rho^2\delta^2+\varepsilon^2}{\delta^2+\varepsilon^2}\leq\rho^2
$$
then, from definition  \eqref{defVstar} of $\mathcal{V^*}$, replacing $\lambda$ and $M$ with $\delta^2+\varepsilon^2$ and $\rho$ respectively, we infer the statement of the proposition.
\qed

Let $V(x)\in \mathcal{V}^*$ and $\bar u(t)$ any collision solution of the system \eqref{sistema} with energy $E$ and leading from the  initial condition $\bar\nu=(\bar q_0,\bar p_0)$ in the form \eqref{iniz1coll} or \eqref{iniz2coll}.
 According with the previous  setting, section \ref{setting},  for every sufficiently small $\varepsilon>0$, let $u_{\varepsilon,\nu}(t)$ be the solution of the regularized system \eqref{sistemareg} with initial condition $\nu$. Reminding the definition of  $\Delta\theta(u_{\varepsilon,\nu})$ given in \eqref{dteta1}, \eqref{dteta2} we state the following theorem.

\begin{teorema}\label{teoregolariz}
There exists
$$
\lim_{(\varepsilon,\nu)\rightarrow (0,\bar \nu)}\Delta\theta(u_{\varepsilon,\nu})
$$
and such limit is $\frac{\pi}{2}$.
\end{teorema}

\dimostr
Reminding the definition of $\bar R$ and $R_+$,  we define
$$
\beta=\min\{\bar R, R_+\}
$$
therefore, for every $u_{\varepsilon,l}(t)$, regardless they are bounded or not by $\bar R$, we write 
$$
\Delta\theta(u_{\varepsilon,\nu})=\int_{R_-}^{\beta}\frac{l}{r^2\dot r}dr
$$
In order to deal with the convergence of $\Delta\theta(u_{\varepsilon,\nu})$ we  
 first rewrite the integrand  in a different way.
From the conservation of energy it follows
$$
\dot r^2=2(E_{\varepsilon,\nu}+V_\varepsilon(r))-\frac{l^2}{r^2}
$$
thus, replacing the radial velocity and extracting the roots $R_-$ and $\beta$ from the denominator, we infer
\begin{align*}
\displaystyle{\Delta}\theta(u_{\varepsilon,\nu})&=\displaystyle{\int_{R_-}^\beta
\frac{1}{r\sqrt{\frac{2r^2}{l^2}(E_{\varepsilon,\nu}+V_\varepsilon(r))-1}}dr}=\displaystyle{\int_{R_-}^\beta\frac{1}{r\sqrt{(r-R_-)(\beta-r)}}
\sqrt{\frac{(r-R_-)(\beta-r)}{\frac{2r^2}{l^2}(E_{\varepsilon,\nu}+V_\varepsilon(r))-1}}}
\end{align*}
By means of the change of variables $\rho=\frac{r}{R_-}$ we rewrite the integral in the form
 \begin{gather}
{\displaystyle \Delta\theta(u_{\varepsilon,\nu})=\int_1^\frac{\beta}{R_-}
\frac{1}{\rho\sqrt{(\beta-\rho R_-)(\rho-1)}}\sqrt{\frac{(\beta-\rho R_-)(\rho-1)}{\frac{2R_-^2\rho^2}{l^2}\Big(E_{\varepsilon,\nu}+V_\varepsilon(\rho R_-)\Big)-1}}d\rho}\nonumber\\
{\displaystyle =\int_1^\frac{\beta}{R_-}
\frac{1}{\rho\sqrt{(\bar R-\rho R_-)(\rho-1)}}\sqrt{K}d\rho}\label{farelimite2}
\end{gather}

We proceed with the proof of the theorem as follows: first we exhibit, for every $V(x)\in\mathcal{V}$, an uniform bound for the function $K$ provided $\varepsilon$ small enough is taken, then we apply the Lebesgue's theorem in order to  obtain the limit of $\Delta\theta(u_{\varepsilon,\nu})$ as
 $(\varepsilon,\nu)\rightarrow
(0,\bar \nu)$ and we'll prove that such limit exists if $V(x)\in\mathcal{V}^*$.

\begin{lemma}\label{limitatezzak1k2}
Let $V(x)\in\mathcal{V}$, then for $\varepsilon$ small enough the function $K$ is bounded by a constant in its domain.
\end{lemma}
\dimostr
Replacing in $K$ the relation $E_{\varepsilon,\nu}=\frac{1}{2}v^2+\frac{1}{2}\frac{l^2}{\beta^2}-V_\varepsilon(\beta)$ 
we obtain 
\begin{equation}
K=\frac{(\beta-\rho R_-)(\rho-1)}{\frac{2R_-^2\rho^2}{l^2}\Big(\frac{1}{2}{v^2+\frac{1}{2}\frac{l^2}{\beta^2}
-V_\varepsilon(\beta)+V_\varepsilon(\rho R_-)\Big)-1}}
\end{equation}
We observe that $v$ is the zero in case $\beta=R_+$.

Subtracting the energy formula $E_{\varepsilon,\nu}=\frac{1}{2}|\dot u_{\varepsilon,l}|^2-V_\varepsilon(u_{\varepsilon,l})$ evaluated in $R_-$ from the same evaluated in $\beta$ we infer
\begin{equation*}
\frac{2R_-^2}{l^2}=\frac{1}{(V_\varepsilon(R_-)-V_\varepsilon(\beta)
+\frac{1}{2}v^2)}\frac{\beta^2-R_-^2}{\beta^2}
\end{equation*}
that, replaced into $K$, implies 
\begin{gather}
\displaystyle{K=\frac{(\beta-\rho R_-)(\rho-1)}{\left(\frac{R_-^2\rho^2}{\beta^2}-1\right)+\rho^2\left(\frac{V_\varepsilon(\rho R_-)-V_\varepsilon(\beta )+\frac{v_1^2}{2}}{V_\varepsilon(R_-)-V_\varepsilon(\beta)+\frac{v_1^2}{2}}\right)\frac{\beta ^2-R_-^2}{\beta^2}}}\label{usareK}
\end{gather}

We note that 
for every $0<a<b$, the function $f(z)=\frac{a+z}{b+z}$ is increasing for positive $z$. The condition \textbf{ii.} in definition \eqref{defV} implies that, for every small enough $\varepsilon$, the function $V_\varepsilon(x)$ is decreasing with respect to $x$ for every $x<\beta$. This yields the relations
$V_\varepsilon(\beta)<V_\varepsilon(\rho R_-)<V_\varepsilon(R_-)$, thus replacing $a=V_\varepsilon(\rho R_-)-V_\varepsilon(\beta)$, $b=V_\varepsilon(R_-)-V_\varepsilon(\beta)$ and $z=\frac{v_1^2}{2}$ in $f(z)$, it follows 
$$
\frac{V_\varepsilon(\rho R_-)-V_\varepsilon(\beta )+\frac{v^2}{2}}{V_\varepsilon(R_-)-V_\varepsilon(\beta)+\frac{v^2}{2}}>\frac{V_\varepsilon(\rho R_-)-V_\varepsilon(\beta)}{V_\varepsilon(R_-)-V_\varepsilon(\beta)}\qquad \forall \rho\in\left(1,\frac{\beta}{R_-}\right),\ \forall v^2
$$
Therefore
\begin{eqnarray}
&K\displaystyle{<\frac{\beta^2}{\frac{(\beta+\rho R_-)}{\rho-1}\left(\rho^2\Big(\frac{V_\varepsilon(\rho R_-)-V_\varepsilon(\beta)}{V_\varepsilon(R_-)-V_\varepsilon(\beta)}\Big)
\frac{\beta^2-R_-^2}{\beta^2-R_-^2\rho^2}-1\right)}}\nonumber
\end{eqnarray}
and, by means of  the  substitutions $x=\rho R_-$ and $y=R_-$

$$
K<\displaystyle{\frac{\beta^2}{\frac{(\bar R+x)}{\frac{x}{y}-1}\left(\left(\frac{x}{y}\right)^2
\Big(\frac{V_\varepsilon(x)-V_\varepsilon(\beta)}
{V_\varepsilon(y)-V_\varepsilon(\beta)}\Big)\frac{\beta^2-y^2}{\beta ^2-x^2}-1\right)}} \qquad x\in(y,\beta)
$$
Finally, applying proposition (\ref{limitatezzaK}), we obtain the uniform bound
$$
K<\beta\quad \forall \rho\in\left(1,\frac{\beta}{R_-}\right)    
$$
\qed

\textit{Continue the proof of theorem \ref{teoregolariz}.}

\noindent
The boundness of the functions $K_1$ and $K_2$ and the formula
\begin{eqnarray}\label{fapig}
\int_1^\xi\frac{1}{x\sqrt{(x-1)(1-\frac{x}{\xi})}}dx=\pi\quad \xi>1
\end{eqnarray}
implies that all the integral  \eqref{farelimite2} is bounded by $\pi\beta$.

In order to apply the Lebesgue dominated convergence theorem we need an uniform bound of all the  integrands, independently on the values of $(\varepsilon,\nu)$ in the neighbourhood of $(0,\bar \nu)$. We note that  the singular point $\frac{\beta}{R_-}$ moves as $\varepsilon$ and $\nu$ change, thus we rewrite   $\Delta\theta(u_{\varepsilon,\nu})$
as follows
\begin{gather*}
\Delta\theta(u_{\varepsilon,\nu})=\displaystyle{\int_1^{\frac{\beta}{R_-}}\frac{1}{\rho\sqrt{(\beta-\rho R_-)(\rho-1)}}\sqrt{K}d\rho}\\
=\displaystyle{\int_1^{\sqrt{\frac{\beta}{R_-}}}\frac{1}{\rho\sqrt{(\beta-\rho R_-)(\rho-1)}}\sqrt{K}d\rho}
+\int_{\sqrt{\frac{\beta}{R_-}}}^\frac{\beta}{R_-}\frac{1}{\rho\sqrt{(\beta-\rho R_-)(\rho-1)}}\sqrt{K}d\rho=I_1+I_2
\end{gather*}
and
\begin{equation*}
\lim_{(\varepsilon,\nu)\rightarrow(0,\bar \nu)}\Delta\theta(u_{\varepsilon,\nu})=
\lim_{(\varepsilon,\nu)\rightarrow(0,\bar \nu)}I_1+\lim_{(\varepsilon,\nu)\rightarrow(0,\bar \nu)}I_2
\end{equation*}
First we calculate the integral $I_2$ and we check that it  is infinitesimal as $({\varepsilon,\nu})\rightarrow(0,\bar \nu)$:
for the boundness of $K$, $K\leq\beta$,
\begin{align*}
I_2&=\displaystyle{\frac{1}{\sqrt{\beta}}
\int_{\sqrt{\frac{\beta}{R_-}}}^\frac{\beta}{R_-}
\frac{1}{\rho\sqrt{(1-\rho\frac{R_-}{\beta})(\rho-1)}}
\sqrt{K}d\rho}<\displaystyle{\int_{\sqrt{\frac{\beta}{R_-}}}
^\frac{\beta}{R_-}\frac{1}{\rho\sqrt{(1-\rho\frac{R_-}{\beta})
(\rho-1)}}d\rho}
\end{align*}
hence, for \eqref{fapig}, 
\begin{align*}
I_2&<\pi-\displaystyle{\int_1^{\sqrt{\frac{\beta}{R_-}}}
\frac{1}{\rho\sqrt{(1-\rho\frac{R_-}{\beta})(\rho-1)}}d\rho}
=\pi-2\arctan\left[\sqrt{\frac{\rho-1}{1-\frac{\rho R_-}{\beta}}}\,\right]\Biggr|^{\sqrt{\frac{\beta}{R_-}}}_1\\
&=2\left(\frac{\pi}{2}-\arctan\sqrt{\frac{\sqrt{\beta}}{\sqrt{R_-}}}\right)=2\arctan\sqrt{\frac{\sqrt{R_-}}{\sqrt{\beta}}}\approx 2\sqrt[4]{\frac{R_-}{\beta}}
\end{align*}
where, in the last passage, the relation $\arctan(x)+\arctan(\frac{1}{x})=\frac{\pi}{2}$ is used.
Passing to the limit, reminding that $R_-\rightarrow 0$ as $({\varepsilon,\nu})\rightarrow(0,\bar \nu)$, we infer
\begin{equation*}
\lim_{(\varepsilon,\nu)\rightarrow(0,\bar \nu)}I_2=0\ .
\end{equation*}
It remains to prove the existence of the limit for $I_1$.
We define  $F_{{\varepsilon,\nu}}(\rho):\mathbb{R}^+\rightarrow\mathbb{R}^+$  the function
$$
F_{\varepsilon,\nu}(\rho):=\frac{1}{\rho\sqrt{(\rho-1)(\beta-\rho R_-)}}\sqrt{K}\ \chi_{\left[1,\sqrt{\frac{\beta}{R_-}}\right]}
$$
thus the limit of $I_1$ is equivalent to
$$
\lim_{(\varepsilon,\nu)\rightarrow(0,\bar \nu)}\int_1^\infty F_{\varepsilon,\nu}(\rho)d\rho
$$
We note that every function $F_{{\varepsilon,\nu}}(\rho)$ becomes  unbounded only for $\rho$ approaching $\rho=1$. 
Moreover, for the boundness of $K$, it follows that $\forall \rho\in(1,+\infty)$
\begin{align*}
F_{\varepsilon,\nu}(\rho)&<\frac{C}{\rho\sqrt{(\rho-1)(1-\frac{\rho R_-}{\beta})}}\chi_{\left[1,\sqrt{\frac{\beta}{R_-}}\right]}
<\frac{C}{\rho\sqrt{(\rho-1)\Big(1-\sqrt{\frac{R_-}{\beta}}\Big)}}\chi_{\left[1,\sqrt{\frac{\beta}{R_-}}\right]}<\frac{C'}{\rho\sqrt{\rho-1}}\quad 
\end{align*}
then all the functions  $F_{\varepsilon,\nu}(\rho)$ are dominated by a function $\tilde{F}\in\mathbb{L}^1([1,\infty])$. We apply the Lebesgue theorem computing the pointwise limit of $F_{\varepsilon,\nu}(\rho)$ as $(\varepsilon,\nu)\rightarrow(0,\bar \nu)$.

Using $K$ in the form \eqref{usareK} we write
 \begin{eqnarray}
&\displaystyle{\lim_{(\varepsilon,\nu)\rightarrow(0,\bar \nu)}\frac{1}
{\rho\sqrt{(\rho-1)(\beta-\rho R_-)}}\sqrt{K_1}\chi_{\left[1,\sqrt{\frac{\beta}{R_-}}\right]}}= \nonumber\\
&\displaystyle{\lim_{(l,\varepsilon)\rightarrow(0,0)}
\frac{1}{\rho}\sqrt{\frac{\beta^2}{(\beta^2-R_-^2\rho^2)\left(\rho^2
\frac{\big(V_\varepsilon(\rho R_-)-V_\varepsilon(\beta)+\frac{v^2}{2}\big)}
{\big(V_\varepsilon(R_-)-V_\varepsilon(\beta)+\frac{v^2}{2}\big)}\frac{(\beta^2-R_-^2)}{(\beta^2-R_-^2\rho^2)}-1\right)}}
\chi_{\left[1,\sqrt{\frac{\beta}{R_-}}\right]}}\nonumber
\end{eqnarray}
and reminding that the convergence $(\varepsilon,\nu)\rightarrow(0,\bar \nu)$ implies $l\rightarrow 0$ and, as a consequence, $R_-\rightarrow 0$ we gain 
$$
F_{\varepsilon,l}(\rho)\sim\frac{1}{\rho}\sqrt{\frac{1}{\rho^2\frac{V_\varepsilon(\rho R_-)}{V_\varepsilon(R_-)}-1}}
$$
Therefore for every $V(x)\in\mathcal{V}^*$, thank to proposition \eqref{allimite}, we infer

$$
\lim_{(l,\varepsilon)\rightarrow(0,0)}F_{\varepsilon,l}(\rho)=\frac{1}{\rho\sqrt{\rho^2-1}}
$$
uniformly in $\varepsilon$ and $l$.
Thus, for the Lebesgue theorem, there exists the limit of $\Delta\theta(u_{\varepsilon,\nu)}$ and it values 
$$
\lim_{(l,\varepsilon)\rightarrow(0,0)}\int_1^\infty F_{\varepsilon,l}(\rho)d\rho=\frac{\pi}{2} .
$$

\qed

Since for $(\varepsilon,\nu)\rightarrow(0,\bar \nu)$ the apsidal angle $\Delta\theta( u_{\varepsilon,\nu})$ tends to $\frac{\pi}{2}$, the pointwise limit of the sequence of trajectories $u_{\varepsilon,l}(t)$ in the ball $B_0(\bar R)$ is a straight line trajectory that crosses the origin.

This suggests to  extend  the collision solution $\bar u(t)$  beyond the singularity replacing symmetrically  the solution itself forward the collision point in the same direction.

\begin{definizione}\label{def-trasm}
 Let  $\bar u(t)$, $t\in [0,T_0)$, be a collision path, $T_0$ the collision instant. Define the \textbf{transmission solution} $u_0(t)$, $t\in[0,2T_0]$ as
\begin{gather*}
\left\{\begin{array}{ll}
u_0(t)=\bar u(t) &t\in[0,T_0]\\
u_0(t)=-\bar u(2T_0-t)&t\in[T_0,2T_0]
\end{array}
\right.
\end{gather*} 
\end{definizione}

 In order to complete the proof of theorem \ref{Mainteo} it remains to show that, for every $t\in[0,2T_0]$, the sequence $\{u_{\varepsilon,\nu}(t)\}$ pointwise  converges to $u_0(t)$ as $({\varepsilon, \nu})\rightarrow(0,\bar \nu)$ and that the flow obtained replacing the collision solution $\bar u$ with the transmission solution $u_0$ is continuous with respect initial data.

\subsection{The regularity of the extended flow}\label{secreg}

We denote with $\mathcal{F}=\{y=(x,\dot x)\in \mathbb{R}^2\times\mathbb{R}^2 \}$  the phase space of planar motion and we consider the initial value problems defined on $\mathcal{F}$ equivalent to systems \eqref{sistema} and \eqref{sistemareg}
\begin{equation*}
P(0)=\left\{\begin{array}{ll}
&{\displaystyle y'=f(y)} \\
& {\displaystyle x(0)=x_0  \in \mathbb{R}^2 \backslash\{ 0\}}\\
& {\displaystyle \dot x(0)=\dot x_0  \in \mathbb{R^2} }
\end{array} \right.
,\ P(\varepsilon)=\left\{\begin{array}{ll}
&{\displaystyle y'=f_\varepsilon(y)} \\
& {\displaystyle x(0)=x_0  \in \mathbb{R^2}} \\
& {\displaystyle \dot x(0)=\dot x_0 \in \mathrm{R^2} }
\end{array} \right.
\end{equation*}
where $f(x,\dot x)=(\dot x,\nabla V(|x|))$ and $f_\varepsilon(x,\dot x)=(\dot x,\nabla V_\varepsilon(|x|))$. 

For every initial data $\bar y=(\bar q,\bar p)$, $|\bar q|\leq\bar R$, leading to collision for the system $P(0)$, let $\bar y(t)=~(\bar x(t),\dot{ \bar x}(t)):[0,T_0)\rightarrow\mathcal{F}$ be the corresponding singular solution where $T_0$ denote the collision time and $|\bar x(t)|\leq\bar R$ for every $t\in[0,T_0)$.
We extend $\bar y(t)$ according to the definition \ref{def-trasm} defining $y_0(t)=(x_0(t),\dot x_0(t))$ as
\begin{equation*}
\begin{array}{ll}
y_0(t)=\bar y(t) &t\in[0,T_0)\\
\\
y_0(t)=\left\{
\begin{array}{l}
x_0(t)=-\bar x(2T_0-t)\\
\dot x_0(t)=\dot{\bar x}(2T_0-t)
\end{array}
\right. &t\in(T_0,2T_0)\\
\end{array} 
\end{equation*}
 The  extension for the collision solutions allows to define the flow beyond the singularity, therefore for $\varepsilon\geq 0$ we denote with  $\Phi_\varepsilon(y,t):\mathcal{F}\times\mathbb{R}^+\rightarrow\mathcal{F}$ the flow associated to the system $P(\varepsilon)$.

Our aim is to study the continuity of $\Phi_\varepsilon(y,t)$ with respect to initial data $y\in\mathcal{F}$ and $\varepsilon$ as $(y,\varepsilon)\rightarrow(\bar y,0)$; to this end we consider the Poincar\'e map  $\Phi_T(y,\varepsilon)$ defined as the solution at time $T$ of the system $P(\varepsilon)$ with initial value $y$ and we show that
$$
\lim_{\substack{
y\rightarrow \bar y\\
\varepsilon\rightarrow 0
}}\Phi_T(y,\varepsilon)=\Phi_T(\bar y,0)
$$
for every $T\not =T_0$.

We note that the continuity of the Poincar\'e map implies the proof of theorem \ref{Mainteo}.
\begin{osservazione}
We can not expect the continuity of the Poincar\'e map in $T_0$ because, even if the configurations $x_{T_0}(y,\varepsilon)$ would converge to $x_{T_0}(\bar y,0)$, the limit can not be attained by the sequence $\dot x_{T_0}(y,\varepsilon)$, since  $\dot x_{T_0}(\bar y,0)$ is unbounded.
\end{osservazione}

If $T<T_0$ no problem arises, indeed the above limit comes from the classical theorem of continuity with respect initial data of ordinary differential equations.   

On the other hand, for $T>T_0$ the continuity of the Poincar\'e map is stated in the following theorem.

\begin{teorema}\label{teopoincare}
Let $V(x)\in\mathcal{V}^*$ and suppose $\bar y=(\bar q,{\bar p})$, $|\bar q|\leq \bar R$  be 
 an initial condition leading to collision for the system $P(0)$ at time $T_0$.
Then
$$
\lim_{\substack{
y\rightarrow \bar y\\
\varepsilon\rightarrow 0
}}\Phi_T(y,\varepsilon)=\Phi_T(\bar y,0)
$$
for every $T>T_0$ such that $\dot x(T)\not =0$
.
\end{teorema}

In the proof of the theorem we will need the following classical lemma
\begin{lemma}\label{implicita}
Let $H(r,p)$ and $H_0(p)$ be real continuous functions with respect to a set of parameters $p$, and suppose $H$ be strictly increasing as function of $r$. Then for every $T$ the function $r=r(T,p)$, implicit solution of equation
$$
T=H_0(p)+H\Big(r(T,p),p\Big)
$$
is continuous as function of $p$.
\end{lemma}

{\it Proof of theorem \eqref{teopoincare}}

\noindent As usual we set $(r,\theta)$ the polar coordinates of the plane  then a point $y=(x,\dot x)$ in the phase space is replaced by  
$$x=(r\cos\theta,r\sin\theta)
$$
$$\dot x=(\dot r\cos \theta-r\dot \theta \sin\theta,\dot r\sin\theta+r\dot \theta\cos\theta)
$$
From the definition of  $\Phi_T(y,\varepsilon)$ remain well defined the functions  and $r_T(y,\varepsilon),\theta_T(y,\varepsilon)$, $\dot r_T(y,\varepsilon),\dot \theta_T(y,\varepsilon)$ denoting, respectively, the values of the radial and angular coordinate and their velocity at time $T\not=T_0$ of a solution of system $P(\varepsilon)$ 
with initial data $y$. The continuity of $\Phi_T(y,\varepsilon)$ is equivalent 
to the continuity of each of the previous functions.

Let $\bar y=(\bar r,\bar \theta,\dot {\bar r},\dot {\bar\theta})$ be an initial data in the phase space  leading to collision with nonzero initial velocity, $\dot{\bar r}<0$, and denote with $\bar E=\frac{1}{2}\dot {\bar r}^2-V(\bar r)$ the energy of the collision  solution.

A point $y\in\mathcal{F}$, $y=(r_0,\theta_0,\dot r_0,\dot\theta_0)$, tends to $\bar y$ if it holds
\begin{gather*}
|r_0-\bar r|\rightarrow 0\ ,\qquad |\theta_0-\bar \theta|\rightarrow 0 \quad \rm{mod}\  2\pi\\
(E,l)\rightarrow(\bar E,0)
\end{gather*}

We start showing the continuity of the function $r_T(y,\varepsilon)$ as $y\rightarrow \bar y$ and $\varepsilon\rightarrow 0$.
The value of $r_T(y,\varepsilon)$ is governed by the differential equation $\dot r=\sqrt{2(E+V_\varepsilon(r))-\frac{l^2}{r^2}}$, thus $r_T(y,\varepsilon)$ is a function of the initial position $r_0$, the couple 
 $E$, $l$ and the parameter $\varepsilon$. Define $\mathcal{T}_0(E,l, r_0,\varepsilon)$ as  the time necessary to the solution $r(y,t)$ to reach the minimal value $R_-=R_-(E,l)$, then
\begin{align*}
T&=\int_{R_-}^{ r_0}\frac{1}{\sqrt{2(E+V_\varepsilon(\rho))-\frac{l^2}{\rho^2}}}d\rho+
\int_{R_-}^{r_T(y,\varepsilon)}\frac{1}{\sqrt{2(E+V_\varepsilon(\rho))-\frac{l^2}{\rho^2}}}d\rho\\
&=\mathcal{T}_0(E,l,r_0,\varepsilon)+\mathcal{T}(r_T(y,\varepsilon),E,l)
\end{align*}

\underline{Claim}
\textsl{The function $\mathcal{T}_0$ and $\mathcal{T}$ are continuous respect to $r_0,E,l,\varepsilon$}.

Suppose for the moment that the claim is true, since
the function $\mathcal{T}$ is strictly increasing with respect to $r_T(y,\varepsilon
)$, for lemma \eqref{implicita}, the function $r_T(y,\varepsilon)$ is continuous with respect to the set of parameters  $E,l,r_0,\varepsilon$. Therefore 
$$ 
\lim_{\substack{
y\rightarrow \bar y\\
\varepsilon\rightarrow 0
}}r_T(y,\varepsilon)=r_T(\bar y)\ .
$$
The continuity of $\theta_T(y,\varepsilon)$ is equivalent to the continuity of $\Delta\theta_T(y,\varepsilon):=\theta_T(y,\varepsilon)-\theta_0$. For the definition of transmission  solution $\Delta\theta_T(\bar y)=\pi$, while, in theorem \eqref{teoregolariz}, we showed that the apsidal angle associated to a solution of the perturbed system tends to $\frac{\pi}{2}$ whenever the initial data tends to a colliding one and the parameter $\varepsilon$ tends to zero. Thus we gain
$$ 
\lim_{\substack{
y\rightarrow \bar y\\
\varepsilon\rightarrow 0
}}\theta_T(y,\varepsilon)=\theta_T(\bar y)
$$

The continuity of $\dot r_T(y,\varepsilon)$ and $\dot\theta_T(y,\varepsilon)$ follows immediately by the continuity of $r_T(y,\varepsilon)$ and relations
\begin{gather*}
\dot r_T(y,\varepsilon)=2\sqrt{E+V_\varepsilon(r_T(y,\varepsilon))-
\frac{1}{2}\frac{l^2}{r_T(y,\varepsilon)^2}}\, , \qquad\dot \theta_T(y,\varepsilon)=\frac{l}{r_T(y,\varepsilon)^2}
\end{gather*} 
\qed

\vspace{0.5cm}
{\it Proof of the claim}

\noindent Denoting with $p$ and $\bar p$ the sets of parameters $p=(E,l,r_0,\varepsilon)$ and $\bar p=(\bar E,0,\bar r,0)$, we have to show that
$$
\lim_{p\rightarrow\bar p}\mathcal{T}_0(p)=\mathcal{T}_0(\bar p)\ .
$$
We want to apply the dominated convergence theorem and pass the limit under the integral sign in
$$
\lim_{p\rightarrow \bar p}\int_{R_-}^{r_0}\frac{1}{\sqrt{2(E+V_\varepsilon(\rho))-\frac{l^2}{\rho^2}}}d\rho
$$ 
To this aim we first exhibit an $\mathbb{L}^1$ bound for the integrand function. We observe that the integrand become  unbounded only for $\rho$ approaching to $R_-$ since we have chosen the initial value $\dot{\bar r}<0$ and, as a consequence, we can suppose $\dot r_0<0$.

Let $\beta\in(R_-,r_0)$ and write
\begin{align}
\int_{R_-}^{r_0}\frac{1}{\sqrt{2(E+V_\varepsilon(\rho))
-\frac{l^2}{\rho^2}}}d\rho&=
\int_{R_-}^{\beta}\frac{1}{\sqrt{2(E+V_\varepsilon(\rho))-\frac{l^2}{\rho^2}}}d\rho+
\int_{\beta}^{ r_0}\frac{1}{\sqrt{2(E+V_\varepsilon(\rho))-\frac{l^2}{\rho^2}}}d\rho\nonumber\\
&=I+II\nonumber\label{farelimite}
\end{align}
Since for every $\rho\in [\beta,r_0]$ the integrand is bounded, it follows
$
II\leq C_1(r_0-\beta)
$.
In order to show the $\mathbb{L}^1$ boundness of the first integrand, using the definition of the energy integral, we   rewrite I in the form
\begin{align*}
I&=\int_{R_-}^\beta\frac{1}{\sqrt{\dot r_0^2+\frac{l^2}{ r_0^2}-2V_\varepsilon(r_0)+2V_\varepsilon(\rho)-\frac{l^2}{\rho^2}}}d\rho=\int_{R_-}^\beta\frac{1}{\sqrt{\dot r_0^2+l^2\left(\frac{r^2- r_0^2}{r^2 r_0^2}\right)+2\big(V_\varepsilon(\rho)-V_\varepsilon( r_0)\big)}}d\rho
\end{align*}
Again from the definition of energy, evaluated in $r_0$ and in the pericentre $R_-$, 
$$
E=\frac{1}{2}\dot r_0^2+\frac{1}{2}\frac{l^2}{ r_0^2}-V_\varepsilon(r_0)=\frac{1}{2}\frac{l^2}{R_-^2}-V_\varepsilon(R_-)
$$
it follows
$$
l^2=\left(\frac{R_-^2-r_0^2}{r_0^2R_-^2}\right)\left(2\big(V(r_0)-V(R_-)\big)-\dot r_ 0^2\right)
$$
and, replacing the last relation in the integrand, we obtain
\begin{align*}
I&=\int_{R_-}^\beta\frac{1}{\sqrt{\dot r_0^2\left[1-\frac{R_-^2( r_0^2-\rho^2)}{\rho^2( r_0^2-R_-^2)}\right]+\frac{R_-^2}{\rho^2}\frac{(r_0^2-\rho^2)}{(r_0^2-R_-^2)}2\big(V_\varepsilon (r_0)-V_\varepsilon(R_-)\big)+2\big(V_\varepsilon(\rho)-V_\varepsilon(r_0)\big)}}d\rho\\
&\leq\int_{R_-}^\beta\frac{1}{\sqrt{2\big(V_\varepsilon(\rho)-V_\varepsilon( r_0)\big)}}\frac{1}{\sqrt{1-\frac{R_-^2}{\rho^2}\frac{( r_0^2-\rho^2)}{( r_0^2-R_-^2)}\frac{V_\varepsilon(R_-)-V_\varepsilon( r_0)}{V_\varepsilon
(\rho)-V_\varepsilon( r_0)}}}d\rho
\end{align*}
For proposition \ref{limitatezzaK} there exists a positive constant $K$ such that, for every $\varepsilon$ small enough,
$$
1-\frac{R_-^2}{\rho^2}\frac{( r_0^2-\rho^2)}{( r_0^2-R_-^2)}\frac{V_\varepsilon(R_-)-V_\varepsilon( r_0)}{V_\varepsilon(\rho)-V_\varepsilon( r_0)}>\frac{K r_0(\rho-R_-)}{r_-( r_0+\rho)+K r_0(\rho-R_-)}
$$
thus
$$
I<C\int_{R_-}^\beta\frac{1}{\sqrt{\rho-R_-}}\frac{1}{\sqrt{2\big(V(\beta)-V(r_0)\big)}}
<C_1\int_{R_-}^\beta\frac{1}{\sqrt{\rho-R_-}}d\rho<\infty
$$
In order to obtain an uniform bound for every choice of $R_-$, we perform the variable change $z=\rho-R_-$ and conclude
$$
I<C_2\int_0^\beta\frac{1}{\sqrt{z}}dz<\infty
$$ 
We now apply the Lebesgue theorem  and we obtain
$$
\lim_{p\rightarrow \bar p}\int_{R_-}^{ r_0}\frac{1}{\sqrt{2(E+V_\varepsilon(\rho))
-\frac{l^2}{\rho^2}}}d\rho=\int_{0}^{\bar r}\frac{1}{\sqrt{2(\bar E+V(\rho))}}d\rho
$$
thus the function $\mathcal{T}_0$ and, for the same reason, the function $\mathcal{T}$ are continuous with respect to the parameter $p$. 

\qed

\begin{figure}[t]
\centering
\includegraphics[scale=0.5]{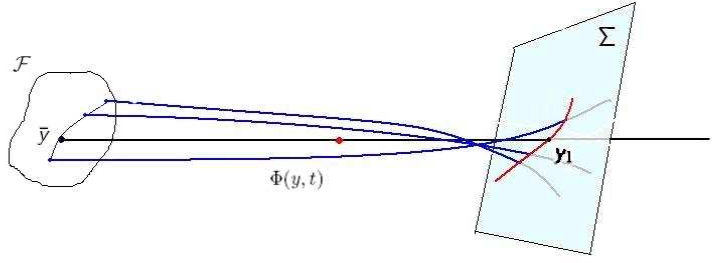}
\caption{Poincar\'e section}
\label{figsezione}
\end{figure}

Let again $\bar y$ be  an initial condition leading to collision for the system $P(0)$, $T_0$ the collision instant and, for a choice of  $T>T_0$, we denote  $y_1=\Phi(\bar y,T)$ according to the definition \ref{def-trasm} of extension of  collision solution.
Given  $\Sigma$ an hyperplane in the phase space passing through $y_1$ and transversally to the flow, we show that for every initial data  $y$ near $\bar y$, there exists a time $t=\tau(y)$ when the trajectory $\Phi(y,t)$ intersects the hyperplane $\Sigma$.  As the data $y$ changes, the trace $S(y)=\Phi(y,\tau(y))$ drawn on $\Sigma$ is called the {\it Poincar\'{e} section} of the flow on $\Sigma$, see figure \ref{figsezione}.
In the next theorem we prove the continuity of the map $\tau(y)$ and the continuity of the Poincar\'{e} section in a neighbourhood of $\bar y $.

\begin{teorema}\label{teocontiuitapoincare}
Let $V(x)\in\mathcal{V}^*$. Then there exists a $\delta>0$ and a continuous function $\tau(y)$ defined in a $\delta$-neighbourhood of $\bar y$, $N_\delta(\bar y)$, such that $\tau(\bar y)=T$ and
$$
\Phi(y,\tau(y))\in\Sigma
$$
Moreover the Poincaré\'{e} section is continuous in $N_\delta(\bar y)$.
\end{teorema}
\dimostr
Let $F$ a vector in the phase space such that $F\cdot f(y_1)>0$ and consider the hyperplane 
$$
\Sigma=\left\{y:(y-y1)\cdot F=0 \right\}
$$
By definition of $y_1$ it holds $(\Phi(\bar y,T)-y1)\cdot F=0$
and, since 
$$
\frac{d}{dt}(\Phi(\bar y,t)-y1)\cdot F\Big|_{t=T}=f(y_1)\cdot F>0
$$
there exists an $\xi>0$ such that
\begin{gather*}
(\Phi(\bar y,T-\xi)-y1)\cdot F<0\\
(\Phi(\bar y,T+\xi)-y1)\cdot F>0
\end{gather*} 
For theorem \eqref{teopoincare} and for the sign permanence theorem it follows that there exists a $\delta>0$ such that for every $|\bar y-y|<\delta$ 
\begin{gather*}
(\Phi(y,T-\xi)-y1)\cdot F<0\\
(\Phi(y,T+\xi)-y1)\cdot F>0
\end{gather*} 
For every fixed $y$ the function $(\Phi(y,t)-y1)\cdot F$ is increasing in $t$: indeed 
$$
\frac{d}{dt}(\Phi(y,t)-y1)\cdot F=f(\Phi(y,t))\cdot F
$$
thus the continuity of the vector space $f(y)$ out of  collision points and the continuity of the orbit $\Phi(y,t)$ with respect to both the variables imply that $f(\Phi(y,t))\cdot F>0$.  
Therefore for every $y\in N_\delta(\bar y)$ there exists a time $\tau(y)$ such that 
$\Phi( y,\tau(y))\in \Sigma$ and $\tau(\bar y)=T$.
In order to prove the continuity of $\tau(y)$ we define 
$$
H(y,\tau)=(\Phi(y,\tau)-y_1)\cdot F
$$
then $\tau(y)$ is the implicit solution of $H(y,\tau(y))=0$.
Since $H$ is continuous in $y$ and it is continuous and increasing with respect to  $\tau$,
 for lemma \eqref{implicita} the function $\tau(y)$ is continuous. Moreover,
for composition of continuous functions we infer the continuity of the Poincar\'{e} section.

\qed

\section{Variational property of the collision solutions}\label{secvar}

In this section  we join a  variational approach that consists in seeking solutions  of the system $\eqref{sistema}$ as critical points of the action functional 
$$
\mathcal{A}(u)=\int_T\mathcal{L}(u,\dot u)dt
$$
where 
$$
\mathcal{L}(u,\dot u)=\frac{1}{2}|\dot u(t)|^2+V(|u(t)|
$$
 is the lagrangian function associated to the equation of motion. This method is well known in the literature and it has been extensively exploited in order to find periodic solutions for the $N$-body problem, see for instance \cite{BFT,FT,Go} and the references therein. Besides  the discussion about the existence of minimal paths, an interesting question is whether the collisions are avoided by the minimal paths even in presence of weak potentials whose contribution one could expect to be negligible by a variational point of view.

In the following theorem we prove that, despite of the weakness  of the singularity of every  potential $V(x)\in\mathcal V^*$, a minimal path for the action functional can not have a collision in the interior of its domain.


\begin{teorema}
For every $V(x)\in\mathcal{V}^*$, let  $u_0:[-T,T]\rightarrow \mathbb{R}^2$ be the extension of a collision  solution of system \eqref{sistema} according with the definition \ref{def-trasm} of transmission solution. Then $u_0(t)$ is not a minimal path for the action functional $\mathcal{A}_{|\Lambda}$, where 
$$
\Lambda=\{u(t):\dot u\in\mathbb{L}^2([-T,T]),u(-T)=u_0(-T),u(T)=u_0(T)\}
$$
denotes the set of paths joining the end points of $u_0$.
\end{teorema}

\dimostr

In order to prove that $u_0(t)$ is not the minima for the action $\mathcal{A}$  among all the paths in $\Lambda$, we perform a variation on the trajectory $u_0$ that removes the collision and makes the action decrease.

Let $u_1(t)=u_0(t)+v^\delta(t)$ where   $v^\delta(t)$ is the standard variation 
\begin{equation*}
v^\delta(t)=\left\{
\begin{array}{cl}
\delta\quad&|t|<T_1\\
\frac{(T-t)}{(T_1-T)}\delta \quad&T_1<|t|<T
\end{array}
\right.
\end{equation*}
directed orthogonally to $u_0(t)$.

Let us compute the difference $\Delta\mathcal{A}=\mathcal{A}(u_0)-\mathcal{A}(u_1)$ and show that for every $\delta$ sufficiently small $\Delta\mathcal{A}$ is positive: this means that $u_0(t)$ is not a local minimum of the action functional.
\begin{gather*}
\Delta\mathcal{A}=\int_{-T}^T\Big(\mathcal{L}(u_0)-\mathcal{L}(u_1)\Big)dt=\int_{-T}^T
\left(\frac{1}{2}|\dot u_0|^2+V(|u_0|)\right)- \left(\frac{1}{2}|\dot u_1|^2+V(|u_1|)\right)dt\\
=\int_{-T}^T\frac{1}{2}\left( |\dot u_0|^2-|\dot u_1|^2\right)dt+
\int_{-T}^T V(|u_0|)-V(|u_1|)dt\\
=\Delta\mathcal{K}+\Delta{V}
\end{gather*}
We study separately the kinetic and the potential contribute.

Since the variation $v^\delta$ is directed orthogonally to $u_0$, we gain 
\begin{equation*}
|\dot u_1(t)|^2=\left\{
\begin{array}{cl}
|\dot u_0(t)|^2 \quad& |t|\leq T_1\\
|\dot u_0(t)|^2+\frac{\delta^2}{(T-T_1)^2}\quad& T_1<|t|\leq T
\end{array}
\right.
\end{equation*}
then
\begin{equation}\label{cinetica}
\Delta\mathcal{K}=-2\int_{T_1}^T\frac{1}{2}\frac{\delta^2}{(T-T_1)^2}dt
=-\frac{\delta^2}{(T-T_1)}
\end{equation}
We show that, for every $\delta$ small enough, the contribution of the potential part is bigger than the penalising contribution, due to the kinetic part. 
\begin{gather*}
\Delta{V}=2\int_0^T V(|u_0|)-V(|u_1|)\,dt>2\int_0^{T_1} V(|u_0|)-V(|u_1|)\,dt\\
\geq \int_0^{T_1}V(|u_0|)-V\big(\sqrt{u_0^2+\delta^2}\big)\,dt 
\end{gather*}
For every $t$ fixed
$$
V(|u_0|)-V\Big(\sqrt{u_0^2+\delta^2}\Big)=-\int_0^1\frac{d}{d\xi}V\left( \sqrt{u_0^2+\xi\delta^2}\right)d\xi
$$
hence
\begin{gather}
\Delta{V}\geq 2\int_0^{T_1}\int_0^1-\frac{V'\big(\sqrt{u_0^2+\xi\delta^2} \big)}{2\sqrt{u_0^2+\xi\delta^2}}\delta^2 \,d\xi\, dt=\delta^2\int_0^{T_1}\int_0^1-\frac{V'\big(\sqrt{u_0^2+\xi\delta^2} \big)}{\sqrt{u_0^2+\xi\delta^2}}\,d\xi\, dt\nonumber\\
=\delta^2\int_0^{T_1}f_\delta(t)\,dt= \delta^2R_\delta\ .\label{potenziale}
\end{gather}
For $\delta$ sufficiently small the functions $f_\delta(t)$ are positive and 
 for  Fatou's Lemma 
\begin{align}
\lim_{\delta\rightarrow 0}
\int_0^{T_1}\int_0^1-\frac{V'\big(\sqrt{u_0^2+\xi\delta^2} \big)}{\sqrt{u_0^2+\xi\delta^2}}\,d\xi\, dt&\geq\int_0^{T_1}\frac{-V'(|u_0(t)|)}{|u_0(t)|}\, dt\nonumber\\
&>C\int_0^{T_1}-V'(|u_0(t)|)\ dt=C\lim_{t\rightarrow 0}\dot r(t)=+\infty\label{relfatou}
\end{align}
The last relation holds since the radial differential equation $\ddot r(t)=-V'(r(t))$ is satisfied, being  the collision solution $u_0$ a radial solution.
Combining \eqref{cinetica}, \eqref{potenziale} and \eqref{relfatou} we conclude
$$
\Delta\mathcal{A}=\Delta\mathcal{K}+\Delta\mathcal{V}=
\delta^2\Big(-\frac{1}{T-T_1}+R_\delta \Big)>0
$$
for every $\delta$ small enough.
\qed


\vskip 10pt
\noindent \texttt{r.castelli3@campus.unimib.it, susanna.terracini@unimib.it}
\vskip 5pt
\noindent Universit\`a di Milano Bicocca,\\
Dipartimento di Matematica e Applicazioni,\\
Via R. Cozzi 53, 20125 Milano, Italy.\\

\end{document}